\documentclass[12pt]{article}
\usepackage{latexsym}
 \usepackage{mathrsfs}
\usepackage{amssymb}
\usepackage{graphicx}
\usepackage{cite}
\usepackage{dsfont}

\newtheorem{Theorem}{Theorem}[part]
\newtheorem{Definition}{Definition}[part]
\newtheorem{Proposition}{Proposition}[part]

\newtheorem{Lemma}{Lemma}[part]
\newtheorem{Corollary}{Corollary}[part]
\newtheorem{Remark}{Remark}[part]
\newtheorem{Example}{Example}[part]

\def \ep{\hbox{ }\hfill$\Box$}

\def\reff#1{{\rm(\ref{#1})}}

\begin{document}
\title{A Theoretical Perspective of Solving Phaseless Compressed Sensing via Its Nonconvex Relaxation\thanks{This work was partially supported by the National Natural
Science Foundation of China (Grant No. 11431002).}
}

\author{
Guowei You\thanks{Department of Mathematics, School of Science, Tianjin University, Tianjin 300072 P.R. China; and Department of Mathematics and Statistics, Henan University of Science and Technology, Luoyang, 471023  P.R. China. Email: youziyou1017@126.com}
\and
Zheng-Hai Huang\thanks{Corresponding Author. Department of Mathematics, School of Science, Tianjin University, Tianjin 300072, P.R. China. This author is also with the Center for Applied Mathematics of Tianjin University. Email: huangzhenghai@tju.edu.cn.   Tel:+86-22-27403615 Fax:+86-22-27403615}
\and
You Wang\thanks{Department of Mathematics, School of Science, Tianjin University, Tianjin 300072, P.R. China. Email: wang\underline{ }yong@tju.edu.cn.}
}

\date{}

\maketitle

\begin{abstract}
\noindent

As a natural extension of compressive sensing and the requirement of some practical problems, Phaseless Compressed Sensing (PCS) has been introduced and studied recently. Many theoretical results have been obtained for PCS with the aid of its convex relaxation. Motivated by successful applications of nonconvex relaxed methods for solving compressive sensing, in this paper, we try to investigate PCS via its nonconvex relaxation. Specifically, we relax PCS in the real context by the corresponding $\ell_p$-minimization with $p\in (0,1)$. We show that there exists a constant $p^\ast\in (0,1]$ such that for any fixed $p\in(0, p^\ast)$, every optimal solution to the $\ell_p$-minimization also solves the concerned problem; and derive an expression of such a constant $p^\ast$ by making use of the known data and the sparsity level of the concerned problem. These provide a theoretical basis for solving this class of problems via the corresponding $\ell_p$-minimization.
\vspace{3mm}

\noindent {\bf Key words:}\hspace{2mm} Phase retrieval, phaseless compressed sensing, compressed sensing, $\mathbf{\ell_p}$-minimization. \vspace{3mm}

\end{abstract}

\section{Introduction}
\setcounter{equation}{0} \setcounter{Assumption}{0}
\setcounter{Theorem}{0} \setcounter{Proposition}{0}
\setcounter{Corollary}{0} \setcounter{Lemma}{0}
\setcounter{Definition}{0} \setcounter{Remark}{0}
\setcounter{Algorithm}{0}

\hspace{4mm} In the past decade,  Compressive Sensing (CS) has gained intensive attention, see \cite{Cades05OnCS,Candes08CS,Donoho05OnCSL1,Herzet13Partial-support,Theodo12OnCS-overview} and references therein. It aims at recovering a sparsest vector from an underdetermined system of linear equations. In other words, CS is to solve the following $\ell_0$-minimization:
\begin{eqnarray}\label{CSoriginal_model}
\min \|\mathbf{x}\|_0\quad \textrm{s.t.}\;\; A\mathbf{x}=\mathbf{b},
\end{eqnarray}
where $A\in \mathbb{R}^{m\times n}$ with $m\ll n$ and $\|\mathbf{x}\|_0$ denotes the number of nonzero elements of $\mathbf{x}\in \mathbb{R}^n$. Unfortunately, problem (\ref{CSoriginal_model}) is NP-hard \cite{Nata95linearL0-NP-hard}. To deal with (\ref{CSoriginal_model}), many methods have been developed, and among which the $\ell_1$-minimization
approach is well-known. To fill the gap between the $\ell_0$-minimization and $\ell_1$-minimization, many authors studied the $\ell_p$-minimization (see, for example, \cite{{Chart07Lp},{Davies09Lp},{Foucart09Lp},{Grib03Union},{Lai11Lp},{peng15equivalence_Lp},{Saab08Lp}, {Sun12Lp},{Wang11Lp},{Xu12Lp},ZHZ-13}):
\begin{eqnarray}\label{lp-minimization-for-linear-model}
\min \|\mathbf{x}\|_p^p\quad \textrm{s.t.}\;\; A\mathbf{x}=\mathbf{b},
\end{eqnarray}
where $p\in (0,1)$ and $\|\mathbf{x}\|_p=(\sum_i|x_i|^p)^{1/p}$ is the Schatten-$p$ quasi-norm of $\mathbf{x}$. It has been shown that it needs fewer measurements with small $p$ for exact recovery via the $\ell_p$-minimization than the $\ell_1$-minimization. Recently, Peng, Yue and Li \cite{peng15equivalence_Lp} showed that there exists a constant $p^\ast>0$ such that for any fixed $p\in (0,p^\ast)$, every optimal solution to the $\ell_p$-minimization also solves the corresponding $\ell_0$-minimization. Such an important property was also studied for some related problems (see, for example, \cite{CX-16,fu11EquivalenceL0Lp,MH15}).

In recent years, Phase Retrieval (PR) has been paid wide attention (\cite{balan09Phaseless,balan06OnPR,candes13MatrixCompletion,candes13PhaseLift}). Mathematically, PR refers to
recovering a vector $\mathbf{x}_0\in\mathbb{C}^{n}$ (or $\mathbb{R}^{n}$) from a set of
phaseless measurements $\{b_j=|\langle\phi_j, \mathbf{x}_0\rangle |,
j=1,2,\ldots, m\}$, where $\phi_j\in \mathbb{C}^{n}$ (or $\mathbb{R}^{n}$) for any
$j=1,2,\ldots,m$. PR has been applied to X-ray imaging, crystallography, electron microscopy and so on. In many applications, the vectors to be recovered are often sparse in
certain basis and, in particular, this occurs in some regimes of
X-ray crystallography \cite{Xu14SPR}. The corresponding model is called Phaseless
Compressed Sensing (PCS) \cite{PhaselessCS2016}.
Recently, PCS in the real context has been studied (\cite{GWX-16,PhaselessCS2016,Xu14SPR}) whose model is given by
\begin{eqnarray*}
(P_0)    \quad\quad\min \|\mathbf{x}\|_0\quad \textrm{s.t.}\;\; |\Phi\mathbf{x}|=\mathbf{b},
\end{eqnarray*}
where $ \Phi=(\phi_1,\ldots,\phi_m)^{\top}\in \mathbb{R}^{m\times n}$ with full row rank and $m\ll n$, $\mathbf{b}\in \mathbb{R}^{m}_+$ and $\mathbf{x}\in \mathbb{R}^{n}$. Hereafter, the symbol $|\cdot|$ denotes the component-wise absolute value of a vector, i.e.,  $|\mathbf{u}|=(|u_1|,\ldots,|u_m|)^{\top}$ for $\mathbf{u}=(u_1,\ldots,u_m)^{\top}\in \mathbb{R}^m$, and the superscript $\top$ represents transposition.

As was done in the case of CS, problem $(P_0)$ can be relaxed by
\begin{eqnarray*}
(P_p)\quad\quad \min \|\mathbf{x}\|_p^p\quad
\textrm{s.t.}\;\;|\Phi\mathbf{x}|=\mathbf{b}\;\; \mbox{\rm where}\; p\in (0,1].
\end{eqnarray*}
When $p=1$, problem $(P_p)$ is called as the $\ell_1$-minimization for PCS. The exact recovery conditions for problem $(P_0)$ by using the $\ell_1$-minimization have been studied (see, for example, \cite{GWX-16,PhaselessCS2016,Xu14SPR}), including Strong Restricted Isometry Property and Null Space Property. However, it seems that problem $(P_p)$ with $p\in (0,1)$ has not been studied so far. An interesting question is {\it whether or not there exists a constant $p^\ast\in (0,1]$ such that for any fixed $p\in (0,p^\ast)$, every optimal solution to problem $(P_p)$ also solves problem $(P_0)$}. In this paper, we answer this question. That is, without any additional assumption, we show that there exists a constant $p^\ast\in (0,1]$ such that every optimal solution to problem $(P_p)$ solves problem $(P_0)$ whenever $p\in (0,p^\ast)$, and derive an expression of $p^\ast$ by making use of matrix $\Phi$, vector $\textbf{b}$ and the sparsity level of the concerned problem.

The rest of the paper is organized as follows. In Section 2, we give some basic concepts and results which will be used in later sections. In Section 3, three subsectopns are included. Specifically, in Subsection 3.1, we discuss the finiteness of optimal solutions of problem $(P_0)$; and furthermore, bound all optimal solutions of problem $(P_0)$ by a box set; in Subsection 3.2, we bound the optimal solution set of problem $(P_p)$ by a box set; and in Subsection 3.3, we give reformulations of problems $(P_0)$ and $(P_p)$, respectively. In Section 4, we show that there exists a constant $p^\ast\in (0,1]$ such that problem $(P_0)$ is equivalent to problem $(P_p)$ whenever $p\in (0,p^\ast)$; and derive an expression of $p^\ast$ by making use of matrix $\Phi$, vector $\textbf{b}$ and the sparsity level of the concerned problem. In Section 5, we give two examples to illustrate our theoretical findings. Some conclusions and comments are given in the last section.

Conventions on some notations in this paper. We denote a matrix by a boldface uppercase letter; a vector by a boldface lowercase letter and a real number by a lowercase letter. All vectors are column vectors. Particularly, the vector of all ones is denoted by $\mathbf{1}$ and the vector of all zeros is denoted by $\mathbf{0}$ whose dimensions are up to the content when it appears. For any positive integer $n$, we denote $[n]:=\{1,2,\ldots, n\}$. For simplicity, we use
$(\mathbf{x},\mathbf{y})$ to denote $(\mathbf{x}^{\top},\mathbf{y}^{\top})^{\top}$;  $\mathbb{R}^n_+$ to denote non-negative $n$-dimensional orthant $\mathbb{R}^n_+:=\{(x_1,\ldots,x_n)^{\top}|x_i\geqslant 0, i\in[n]\}$; $sign(\cdot)$ to denote the sign function, i.e., for any scalar $a\geqslant0$, $sign(a)=1$ if $a>0$, and $sign(a)=0$ if
$a=0$. For any vector $\mathbf{u}\geqslant 0$, we use $\mathbf{sign}(\mathbf{u})$ to denote a vector whose $i$-th element being $sign(u_i)$. For two $n$-dimensional vectors $\mathbf{u}$ and $\mathbf{v}$, $\mathbf{u}\circ\mathbf{v}:=(u_1v_1,\ldots,u_nv_n)^{\top}$. For a matrix $\Phi\in\mathbb{R}^{m\times n}$ and any given index set $I\subset [n]$, $I^c:=[n]\backslash I$ is the complement set of $I$, $\#(I)$ denotes the cardinality of $I$, i.e., the number of elements of $I$; $\Phi_I$ and $\mathbf{x}_I$ denote a sub-matrix and a $\#(I)$-dimensional vector constructed by columns of $\Phi$ and elements of $\mathbf{x}$ corresponding to indices of index set $I$, respectively. Define box set $B_\infty(r):=\{\mathbf{x}\in \mathbb{R}^n \mid \|\mathbf{x}\|_{\infty}\leqslant r\}$ for any given $r>0$.

\section{Preliminaries}
\setcounter{equation}{0} \setcounter{Assumption}{0}
\setcounter{Theorem}{0} \setcounter{Proposition}{0}
\setcounter{Corollary}{0} \setcounter{Lemma}{0}
\setcounter{Definition}{0} \setcounter{Remark}{0}
\setcounter{Algorithm}{0}

\hspace{4mm} In this section, we give some basic concepts and derive several simple results, which will be used in later sections
\begin{Definition}\label{def-nondecreasing-vector-valued function}
\begin{itemize}
\item[(i)] Define a function $\mathbf{f}: \mathbb{R}^n_{+}\rightarrow\mathbb{R}^n_+$ by
\begin{eqnarray}\label{func-add-0}
\mathbf{f}(\mathbf{x})=(f_1(x_1),f_2(x_2),\ldots, f_n(x_n))^{\top},
\end{eqnarray}
where $f_i:\mathbb{R}_+\rightarrow \mathbb{R}_+$ for all $i\in[n]$. The (vector-valued) function $\mathbf{f}$ is said to be monotonically nondecreasing if for any two nonnegative vectors $\mathbf{u}$ and $\mathbf{v}$, $\mathbf{u}\leqslant\mathbf{v}$ implies $\mathbf{f}(\mathbf{u})\leqslant\mathbf{f}(\mathbf{v})$.  In other words, if $0\leqslant u_i\leqslant v_i$ implies $f_i(0)\leqslant f_i(u_i)\leqslant f_i(v_i)$ for any $i\in [n]$.
\item[(ii)] Let $F: \mathbb{R}^n_{+}\rightarrow\mathbb{R}_+$. The (real-valued) function $F$ is said to be monotonically nondecreasing if for any two nonnegative vectors $\mathbf{u}$ and $\mathbf{v}$, $\mathbf{u}\leqslant\mathbf{v}$ implies $F(\mathbf{u})\leqslant F(\mathbf{v})$.
\end{itemize}
\end{Definition}

Two monotonically nondecreasing functions are given in the following examples, which will be frequently used in our subsequent analysis.

\begin{Example}\label{on-sign-fun}
Suppose that $\mathbf{f}: \mathbb{R}^n_{+}\rightarrow\mathbb{R}^n_+$ is defined by (\ref{func-add-0}) with $f_i(x_i)=\textrm{sign}(x_i)$ for any $i\in[n]$. We denote this function by $\mathbf{sign}(\mathbf{\cdot})$. It is easy to see that
$\mathbf{sign}(\mathbf{\cdot})$ is a monotonically nondecreasing function.
\end{Example}

\begin{Example}\label{on-lp-fun}
Suppose that $\mathbf{f}: \mathbb{R}^n_{+}\rightarrow\mathbb{R}^n_+$ is defined by (\ref{func-add-0}) with $f_i(x_i)=x_i^p$ for any $i\in[n]$ and $p\in (0,1]$. It is easy to see that this function is monotonically nondecreasing.
\end{Example}

For any given vectors $\mathbf{u}$ and $\mathbf{v}$ with $\mathbf{0}\leqslant\mathbf{u}\leqslant \mathbf{v}$ and a monotonically nondecreasing function $\mathbf{f}$, we have
$$
\mathbf{1}^{\top}\mathbf{f}(\mathbf{u})=\sum\limits_{i=1}^{n}{f_{i}(u_i)}\leqslant
\sum\limits_{i=1}^{n}{f_{i}(v_i)}=\mathbf{1}^{\top}\mathbf{f}(\mathbf{v}).
$$
Therefore, $\mathbf{1}^{\top}\mathbf{f}$ is monotonically nondecreasing.

Let $\mathbb{D}_{R^n}\subseteq \mathbb{R}^n$ be an arbitrarily given non-empty set, and
$$
\tilde{\mathbb{D}}_{R^{2n}}:=\{(\mathbf{x},\mathbf{y})\mid \mathbf{x}\in \mathbb{D}_{R^n}, |\mathbf{x}|\leqslant\mathbf{y}\}.
$$
Considering the following two  optimization problems with non-empty solution sets:
\begin{eqnarray}
(P_{x}):\;\;\min\limits_{\mathbf{x}\in \mathbb{D}_{R^n}}\mathbf{1}^{T}\mathbf{f}(|\mathbf{x}|)\quad
&\mathrm{and}&\quad
(P_{xy}):\;\;\min\limits_{(\mathbf{x},\mathbf{y})\in \tilde{\mathbb{D}}_{R^{2n}}}\mathbf{1}^{T}\mathbf{f}(\mathbf{y}),\nonumber
\end{eqnarray}
we have the following result.
\begin{Lemma}\label{add-thm1-1}
Suppose that $\mathbf{f}: \mathbb{R}^n_{+}\rightarrow \mathbb{R}^{n}_{+}$ is monotonically nondecreasing. Then, problem $(P_{x})$ is equivalent to problem $(P_{xy})$ in the sense that if $x^*$ is an optimal solution of problem $(P_{x})$, then there exists a vector $y^*$ such that $(x^*,y^*)$ solves $(P_{xy})$; and if $(x^*,y^*)$ is an optimal solution of problem $(P_{xy})$, then $x^*$ solves $(P_{x})$.
\end{Lemma}

\noindent {\bf Proof.} Assume that $\mathbf{x}^\ast\in\mathbb{D}_{R^{n}}$ is an optimal solution of problem $(P_{x})$, then
$$
\mathbf{1}^{\top}\mathbf{f}(|\mathbf{x}^\ast|)\leqslant\mathbf{1}^{\top}\mathbf{f}(|\mathbf{x}|),
\quad \forall \mathbf{x}\in \mathbb{D}_{R^{n}}.
$$
Let $\mathbf{y}^\ast=|\mathbf{x}^\ast|$, then $(\mathbf{x}^\ast, \mathbf{y}^\ast)=(\mathbf{x}^\ast, |\mathbf{x}^\ast|)\in\tilde{\mathbb{D}}_{R^{2n}}$. Since $\mathbf{f}$ is monotonically nondecreasing, it follows that
$$
\mathbf{1}^{\top}\mathbf{f}(\mathbf{y}^\ast)=\mathbf{1}^{\top}\mathbf{f}(|\mathbf{x}^\ast|)\leqslant \mathbf{1}^{\top}\mathbf{f}(|\mathbf{x}|)\leqslant\mathbf{1}^{\top}\mathbf{f}(\mathbf{y})
$$
for any $(\mathbf{x},\mathbf{y})\in \tilde{\mathbb{D}}_{R^{2n}}$. That is, $(\mathbf{x}^\ast, \mathbf{y}^\ast)=(\mathbf{x}^\ast, |\mathbf{x}^\ast|)\in\tilde{\mathbb{D}}_{R^{2n}}$ is an
optimal solution of problem $(P_{xy})$.

Conversely, assume that $(\mathbf{x}^\ast,\mathbf{y}^\ast)$ is an optimal solution of problem $(P_{xy})$, then $\mathbf{x}^\ast\in\mathbb{D}_{R^{n}},|\mathbf{x}^\ast|\leqslant\mathbf{y}^\ast$ and
$$
\mathbf{1}^{\top}\mathbf{f}(\mathbf{y}^\ast)\leqslant\mathbf{1}^{\top}\mathbf{f}(\mathbf{y}), \quad \forall (\mathbf{x},\mathbf{y})\in \tilde{\mathbb{D}}_{R^{2n}}.
$$
Furthermore, for any $\mathbf{x}\in \mathbb{D}_{R^{n}}$, let $\mathbf{y}=|\mathbf{x}|$, then $(\mathbf{x},\mathbf{y})\in \tilde{\mathbb{D}}_{R^{2n}}$. Since $\mathbf{f}$ is monotonically nondecreasing, it follows that
$$
\mathbf{1}^{\top}\mathbf{f}(|\mathbf{x}^\ast|)\leqslant \mathbf{1}^{\top}\mathbf{f}(\mathbf{y}^\ast)\leqslant\mathbf{1}^{\top}\mathbf{f}(\mathbf{y}) =\mathbf{1}^{\top}\mathbf{f}(|\mathbf{x}|), \quad \forall \mathbf{x}\in \mathbb{D}_{R^{n}}.
$$
That is, $\mathbf{x}^\ast\in\mathbb{D}_{R^{n}}$ is an optimal solution of problem $(P_x)$.
\ep

\begin{Corollary}\label{add-cro1-1}
Suppose that $\hat{\mathbb{D}}_{R^n}$ is a non-empty subset of $\mathbb{R}^n$, $\mathbf{l}$ and $\mathbf{u}$ satisfying $\mathbf{l}<\mathbf{0}<\mathbf{u}$ are two given vectors in $\mathbb{R}^n$ and $\mathbf{v}=\max\{-\mathbf{l},\mathbf{u}\}$ whose $i$-th component is $\max\{-l_i,u_i\}$. Denote
\begin{eqnarray*}
\begin{array}{l}
\mathbb{D}_{R^{n}}:=\{\mathbf{x}\mid\mathbf{x}\in \hat{\mathbb{D}}_{R^n},\mathbf{l}\leqslant \mathbf{x}\leqslant \mathbf{u}\},\\
\tilde{\mathbb{D}}_{R^{2n}}:=\{(\mathbf{x}, \mathbf{y})\mid\mathbf{x}\in \mathbb{D}_{R^{n}}, |\mathbf{x}|\leqslant \mathbf{y}\},\\
\mathbb{D}_{R^{2n}}:=\{(\mathbf{x}, \mathbf{y})\mid\mathbf{x}\in \mathbb{D}_{R^{n}}, |\mathbf{x}|\leqslant \mathbf{y}, \mathbf{0}\leqslant \mathbf{y} \leqslant \mathbf{v}\}.
\end{array}
\end{eqnarray*}
If $\mathbf{f}: \mathbb{R}_+^{n}\rightarrow \mathbb{R}^{n}_{+}$  is monotonically nondecreasing, then problem
\begin{eqnarray}\label{pro4-1}
\min_{\mathbf{x}\in \mathbb{D}_{R^{n}}}\mathbf{1}^{\top}\mathbf{f}(|\mathbf{x}|)
\end{eqnarray}
is equivalent to problem
\begin{eqnarray}\label{pro4-2}
\min_{(\mathbf{x},\mathbf{y}) \in \mathbb{D}_{R^{2n}}}\mathbf{1}^{\top}\mathbf{f}(\mathbf{y})
\end{eqnarray}
in the sense that if $x^*$ is an optimal solution of problem (\ref{pro4-1}), then there exists a vector $y^*$ such that $(x^*,y^*)$ solves problem (\ref{pro4-2}); and if $(x^*,y^*)$ is an optimal solution of problem (\ref{pro4-2}), then $x^*$ solves problem (\ref{pro4-1}).
\end{Corollary}

\noindent {\bf Proof.}
By Lemma \ref{add-thm1-1},  problem \reff{pro4-1} is equivalent to problem
\begin{equation}\label{pro4-3}
\min_{(\mathbf{x},\mathbf{y}) \in \tilde{\mathbb{D}}_{R^{2n}}}\mathbf{1}^{\top}\mathbf{f}(\mathbf{y}).
\end{equation}
So, it is sufficient to show that problem (\ref{pro4-3}) is equivalent
to problem (\ref{pro4-2}).

Assume that $(\mathbf{x}^{\ast},\mathbf{y}^{\ast})$ is an optimal solution of problem (\ref{pro4-3}), then $\mathbf{x}^{\ast}\in \mathbb{D}_{R^{n}}$ and $|\mathbf{x}^{\ast}|\leqslant\mathbf{y}^{\ast}$. Let $\tilde{\mathbf{y}}^{\ast}=|\mathbf{x}^{\ast}|$, then $(\mathbf{x}^{\ast}, \tilde{\mathbf{y}}^{\ast})$ is an optimal solution of problem (\ref{pro4-2}). Otherwise, suppose that there exists $(\hat{\mathbf{x}},\hat{\mathbf{y}})\in\mathbb{D}_{R^{2n}}\subseteq\tilde{\mathbb{D}}_{R^{2n}}$ such that $\mathbf{1}^{\top}\mathbf{f}(\hat{\mathbf{y}})<\mathbf{1}^{\top}\mathbf{f}(\tilde{\mathbf{y}}^\ast)$. Since $|\hat{\mathbf{x}}|\leqslant\hat{\mathbf{y}}$ and $\mathbf{f}$ is monotonically nondecreasing, it follows that
$$
\mathbf{1}^{\top}\mathbf{f}(|\hat{\mathbf{x}}|) \leqslant\mathbf{1}^{\top}\mathbf{f}(\hat{\mathbf{y}}) <\mathbf{1}^{\top}\mathbf{f}(\tilde{\mathbf{y}}^\ast)
 =\mathbf{1}^{\top}\mathbf{f}(|\mathbf{x}^\ast|)\leqslant \mathbf{1}^{\top}\mathbf{f}(\mathbf{y}^{\ast}),
$$
which contradicts the assumption that $(\mathbf{x}^{\ast},\mathbf{y}^{\ast})$ is an optimal solution of problem (\ref{pro4-3}).

Conversely, assume that $(\mathbf{x}^\ast,\mathbf{y}^\ast)$ is an optimal solution of problem (\ref{pro4-2}). Let $\tilde{\mathbf{y}}^\ast=|\mathbf{x}^\ast|$, then $(\mathbf{x}^\ast, \tilde{\mathbf{y}}^\ast)\in\tilde{\mathbb{D}}_{R^{2n}}$ must be an optimal solution of problem (\ref{pro4-3}). Otherwise, suppose that there exists $(\hat{\mathbf{x}},\hat{\mathbf{y}})\in\tilde{\mathbb{D}}_{R^{2n}}$ such that
$\mathbf{1}^{\top}\mathbf{f}(\hat{\mathbf{y}})<\mathbf{1}^{\top}\mathbf{f}(\tilde{\mathbf{y}}^\ast)$.
Since $ |\hat{\mathbf{x}}|\leqslant\hat{\mathbf{y}}$ and
$\mathbf{f}$ is monotonically nondecreasing, we have
$$
\mathbf{1}^{\top}\mathbf{f}(|\hat{\mathbf{x}}|)
\leqslant\mathbf{1}^{\top}\mathbf{f}(\hat{\mathbf{y}})<
\mathbf{1}^{\top}\mathbf{f}(\tilde{\mathbf{y}}^{\ast}) =\mathbf{1}^{\top}\mathbf{f}(|\mathbf{x}^\ast|) \leqslant\mathbf{1}^{\top}\mathbf{f}(\mathbf{y}^{\ast}).
$$
Let $\hat{\mathbf{\mathbf{z}}}=|\hat{\mathbf{x}}|$, then it is easy to see that $(\hat{\mathbf{x}}, \hat{\mathbf{z}})\in \mathbb{D}_{R^{2n}}$ and
$\mathbf{1}^{\top}\mathbf{f}(\hat{\mathbf{z}})<\mathbf{1}^{\top}\mathbf{f}(\mathbf{y}^\ast)$,
which contradicts the assumption that $(\mathbf{x}^{\ast},\mathbf{y}^{\ast})$ is an optimal solution of problem (\ref{pro4-2}).
\ep

The following two results can be obtained easily; and hence, we omit their proofs.

\begin{Lemma}\label{add-lem2-1}
Suppose that $\mathbf{x}^*\in\mathbb{D}\subseteq \tilde{\mathbb{D}}\subseteq \mathbb{R}^n$ is an optimal solution of $\min\limits_{\mathbf{x}\in \tilde{\mathbb{D}}}f(\mathbf{x})$. Then, $\mathbf{x}^*$ is an optimal solution of $\min\limits_{\mathbf{x}\in \mathbb{D}}f(\mathbf{x})$.
\end{Lemma}

\begin{Lemma}\label{contains-all-solutions}
Suppose that $\mathbb{D}\subseteq \tilde{\mathbb{D}}\subseteq \mathbb{R}^n$. If $\mathbb{D}$ contains all optimal solutions of $\min\limits_{\mathbf{x}\in \tilde{\mathbb{D}}}f(\mathbf{x})$, then $\min\limits_{\mathbf{x}\in \mathbb{D}}f(\mathbf{x})$ is equivalent to $\min\limits_{\mathbf{x}\in \tilde{\mathbb{D}}}f(\mathbf{x})$.
\end{Lemma}

\section{Reformulations of Problems $(P_0)$ and $(P_p)$}
\setcounter{equation}{0} \setcounter{Assumption}{0}
\setcounter{Theorem}{0} \setcounter{Proposition}{0}
\setcounter{Corollary}{0} \setcounter{Lemma}{0}
\setcounter{Definition}{0} \setcounter{Remark}{0}
\setcounter{Algorithm}{0}  \setcounter{Example}{0}

\hspace{4mm} To achieve our main results, we need to confine the solution sets of problems $(P_0)$ and $(P_p)$ in a same bounded set. To this end, we need to discuss the boundedness of the solution sets of problem $(P_0)$ and problem $(P_p)$, respectively.

\subsection{Boundedness of the Solution Set of Problem $(P_0)$}

\hspace{4mm} Consider the classical $\ell_0$-minimization in the case of CS:
\begin{eqnarray}\label{CS_P0}
\min \|\mathbf{x}\|_0\quad \textrm{s.t.}\;\; \Phi\mathbf{x}=\mathbf{b}_\epsilon,
\end{eqnarray}
where $\mathbf{b}_\epsilon\in\mathbb{R}^m$
with some fixed $\epsilon=(\epsilon_1,\ldots,\epsilon_m)^{\top}\in\{-1,1\}^m$. For
problem (\ref{CS_P0}), we assume that every optimal solution has exactly $s$
nonzero components. We denote the solution set of problem (\ref{CS_P0}) by $S_\epsilon$.

\begin{Lemma}\label{linearly_independent_columns}
Given an arbitrarily optimal solution $\hat{\mathbf{x}}\in S_\epsilon$, we denote its support set by $I_{\hat{\mathbf{x}}}$, i.e., $I_{\hat{\mathbf{x}}}=\{i\mid \hat{x}_i\neq0,i\in[n]\}$ with $\#(I_{\hat{\mathbf{x}}})=s$. Then the corresponding sub-matrix $\Phi_{I_{\hat{\mathbf{x}}}}$ has full column rank.
\end{Lemma}

\noindent {\bf Proof.}
Hereafter, for given $\mathbf{x}_I$,  we use $\mathbf{x}\in \mathbb{R}^n$ to denote an expanded $n$-dimensional vector defined by
$$
\mathbf{x}_I:=\mathbf{x}_I\quad\mbox{\rm and}\quad \mathbf{x}_{I^c}:=\mathbf{0}.
$$
Since $\hat{\mathbf{x}}\in S_\epsilon$, it holds that $\mathbf{b}_\epsilon=\Phi_{I_{\hat{\mathbf{x}}}}\hat{\mathbf{x}}_{I_{\hat{\mathbf{x}}}}$.
Suppose that $\Phi_{I_{\hat{\mathbf{x}}}}$ is not full column rank, then there are more than one solution for the linear equation $\Phi_{I_{\hat{\mathbf{x}}}}\mathbf{x}_{I_{\hat{\mathbf{x}}}}=\mathbf{b}_\epsilon$. Thus, there exists $\mathbf{x}_{I_{\hat{\mathbf{x}}}}^\ast$ with
$\mathbf{x}_{I_{\hat{\mathbf{x}}}}^\ast\neq\hat{\mathbf{x}}_{I_{\hat{\mathbf{x}}}}$
such that $\Phi_{I_{\hat{\mathbf{x}}}}\mathbf{x}^\ast_{I_{\hat{\mathbf{x}}}}=\mathbf{b}_\epsilon$.
Let
$$
\mathbf{x}_{I_{\hat{\mathbf{x}}}}^t =t\hat{\mathbf{x}}_{I_{\hat{\mathbf{x}}}}+(1-t)\mathbf{x}_{I_{\hat{\mathbf{x}}}}^\ast,\quad \forall t\in\mathbb{R},
$$
then, $\mathbf{x}_{I_{\hat{\mathbf{x}}}}^t$ is a solution of
$\Phi_{I_{\hat{\mathbf{x}}}}\mathbf{x}_{I_{\hat{\mathbf{x}}}}^t=\mathbf{b}_\epsilon$ for any $t\in\mathbb{R}$.

Since $\mathbf{x}_{I_{\hat{\mathbf{x}}}}^\ast\neq\hat{\mathbf{x}}_{I_{\hat{\mathbf{x}}}}$,
there exists some $i\in I_{\hat{\mathbf{x}}}$ such that
$$
(\mathbf{x}_{I_{\hat{\mathbf{x}}}}^\ast)_i=d\neq(\hat{\mathbf{x}}_{I_{\hat{\mathbf{x}}}})_i=c\neq 0.
$$
Let $t_0:=\frac{d}{d-c}$ and $\mathbf{x}_{I_{\hat{\mathbf{x}}}}^0:=\mathbf{x}_{I_{\hat{\mathbf{x}}}}^{t_0}$, then we have
$$
(\mathbf{x}^0_{I_{\hat{\mathbf{x}}}})_i =(t_0\hat{\mathbf{x}}_{I_{\hat{\mathbf{x}}}}+(1-t_0)\mathbf{x}_{I_{\hat{\mathbf{x}}}}^\ast)_i=0.
$$
Notice that $\hat{\mathbf{x}}$ and $\mathbf{x}^\ast$ have the same support set $I_{\hat{\mathbf{x}}}$. However, the support set of $\mathbf{x}^0$ is a subset of $I_{\hat{\mathbf{x}}}\setminus \{i\}$, which implies that $\mathbf{x}^0$ is a sparser solution of
problem (\ref{CS_P0}) than $\hat{\mathbf{x}}$. Thus, it yields a contradiction to the assumption that there are exactly $s$ nonzero elements for the optimal solution of problem (\ref{CS_P0}). Therefore,
$\Phi_{I_{\hat{\mathbf{x}}}}$ must have full column rank. We complete
the proof.
\ep

\begin{Remark}\label{uniqueness-on-I}
From the above proof, it follows that there do not exist $\hat{\mathbf{x}},\bar{\mathbf{x}}\in S_\epsilon$ with $\hat{\mathbf{x}}\neq \bar{\mathbf{x}}$ such that $I_{\hat{\mathbf{x}}}=I_{\bar{\mathbf{x}}}$.
\end{Remark}
\begin{Corollary}\label{finite-S-epsilon}
Problem {\rm(\ref{CS_P0})} has finitely many optimal solutions.
\end{Corollary}

\noindent {\bf Proof.}
Suppose that each optimal solution of problem (\ref{CS_P0}), say $\hat{\mathbf{x}}$, has $s\; (\leqslant m)$ nonzero components. By Lemma \ref{linearly_independent_columns}, $\Phi_{I_{\hat{\mathbf{x}}}}$ has full column rank. Such sub-matrices of $\Phi$ have at most $C_n^s$. Thus, by Remark \ref{uniqueness-on-I}, it follows that the number of optimal solutions of problem (\ref{CS_P0}) with exactly $s$ nonzero components is no more than $C_n^s$.
\ep

\begin{Corollary}\label{l0-norm-solution-number-is-finite}
Problem $(P_0)$ has finitely many optimal solutions.
\end{Corollary}

\noindent {\bf Proof.}
Assume that the solution set of problem $(P_0)$ is denoted by $S$. It is easy to see that
$$
S\subseteq \bigcup\limits_{\epsilon\in\{-1,1\}^m} S_\epsilon.
$$
Since the number of elements of every $S_\epsilon$ is finite by Corollary \ref{finite-S-epsilon}, it follows from the above formula that the number of elements of $S$ is finite.
\ep

\begin{Theorem}\label{CS_P0_implicitly-bounds}
Suppose that every optimal solution of problem {\rm(\ref{CS_P0})} has exactly $s\;(\leqslant m)$
nonzero components. Then, problem {\rm(\ref{CS_P0})} is equivalent to
\begin{eqnarray}\label{CS_P0-contain-bounds-implicitly}
\min \|\mathbf{x}\|_0\quad \textrm{s.t.}\;\; \Phi\mathbf{x}=\mathbf{b}_\epsilon, \|\mathbf{x}\|_{\infty}\leqslant c^\epsilon_0,
\end{eqnarray}
where
\begin{eqnarray}\label{define-implicit-bounds}
c^\epsilon_0:=\max\limits_{I\subset[n],\#(I)=s}\max\limits_{i\in I}|((\Phi^{T}_{I}\Phi_{I})^{-1}\Phi^{T}_{I}\mathbf{b}_\epsilon)_i|>0
\end{eqnarray}
with every $\Phi_{I}$ being full column rank.
\end{Theorem}

\noindent {\bf Proof.}
Suppose that $\hat{\mathbf{x}}$ is an arbitrarily optimal solution of problem (\ref{CS_P0}) with support set $I_{\hat{\mathbf{x}}}$ and $\#(I_{\hat{\mathbf{x}}})=s$. Then we have that
$\Phi_{I_{\hat{\mathbf{x}}}}\hat{\mathbf{x}}_{I_{\hat{\mathbf{x}}}}=\mathbf{b}_\epsilon$.
By Lemma \ref{linearly_independent_columns},
$\Phi_{I_{\hat{\mathbf{x}}}}$ has full column rank, and hence,
$\hat{\mathbf{x}}_{I_{\hat{\mathbf{x}}}}=(\Phi^{T}_{I_{\hat{\mathbf{x}}}}\Phi_{I_{\hat{\mathbf{x}}}})^{-1}\Phi^{T}_{I_{\hat{\mathbf{x}}}}\mathbf{b}_\epsilon$.
Obviously,
$$\|\hat{\mathbf{x}}\|_\infty=\|\hat{\mathbf{x}}_{I_{\hat{\mathbf{x}}}}\|_\infty
=\max\limits_{i\in
I_{\hat{\mathbf{x}}}}|((\Phi^{T}_{I_{\hat{\mathbf{x}}}}\Phi_{I_{\hat{\mathbf{x}}}})^{-1}\Phi^{T}_{I_{\hat{\mathbf{x}}}}\mathbf{b}_\epsilon)_i|\leqslant
c^\epsilon_0,$$ where $c^\epsilon_0$ is defined by
(\ref{define-implicit-bounds}). Therefore, $\hat{\mathbf{x}}$ is
contained in the box $B_\infty(c^\epsilon_0)$. Furthermore, the
feasible set of problem (\ref{CS_P0-contain-bounds-implicitly}) contains all
optimal solutions of problem (\ref{CS_P0}). By Lemma
\ref{contains-all-solutions}, problem (\ref{CS_P0}) is equivalent to problem
(\ref{CS_P0-contain-bounds-implicitly}). The proof is complete.
\ep

\begin{Theorem}\label{PCS-implicitly-has-bounds}
Suppose that every optimal solution of problem $(P_0)$  has exactly $s\;(\leqslant m)$ nonzero components. Then, problem $(P_0)$ is equivalent to
\begin{eqnarray}\label{PCS_P0-contain-bounds-implicitly}
\min \|\mathbf{x}\|_0\quad \textrm{s.t.}\;|\Phi\mathbf{x}|=\mathbf{b}, \|\mathbf{x}\|_{\infty}\leqslant c_0:=\max\limits_{\epsilon\in\{-1,1\}^m}c^\epsilon_0>0,
\end{eqnarray}
where $c^\epsilon_0$ is defined by (\ref{define-implicit-bounds}).
\end{Theorem}

{\bf Proof,}
Let $S$  be the optimal solution set of problem $(P_0)$, then
$$
S\subseteq \bigcup\limits_{\epsilon\in\{-1,1\}^m} S_\epsilon.
$$
By Theorem \ref{CS_P0_implicitly-bounds}, every $S_\epsilon$ is contained in the box $$
B_\infty(c^\epsilon_0):=\{\mathbf{x}\in \mathbb{R}^n \mid \|\mathbf{x}\|_{\infty}\leqslant c^\epsilon_0\}.
$$
Then, $S$ is contained in $\bigcup_\epsilon B_\infty(c^\epsilon_0)\subseteq B_\infty(c_0)$ with $c_0$ being given by (\ref{PCS_P0-contain-bounds-implicitly}). Hence, the feasible set of problem (\ref{PCS_P0-contain-bounds-implicitly}) contains all optimal solutions of problem $(P_0)$. By Lemma \ref{contains-all-solutions}, problem (\ref{PCS_P0-contain-bounds-implicitly}) is equivalent to problem $(P_0)$. The proof is complete.
\ep

\subsection{Boundedness of Solution Set to Problem $(P_p)$}

\hspace{4mm} We will bound all optimal solutions of problem ($P_p$) via considering the following  problem with any given $p\in (0,1)$:
\begin{eqnarray}\label{CS_Pp}
\min \|\mathbf{x}\|_p^p\quad \textrm{s.t.}\;\; \Phi\mathbf{x}=\mathbf{b}_\epsilon,
\end{eqnarray}
where $\Phi\in\mathbb{R}^{m\times n}$ is full row rank and $\mathbf{b}_\epsilon=\mathbf{b}\circ\epsilon\in\mathbb{R}^m$. We denote the solution set of problem (\ref{CS_Pp}) by $S^p_\epsilon$.

\begin{Lemma}\label{CS-Lp-bounded}
The solution set $S^p_\epsilon$ of problem {\rm(\ref{CS_Pp})} is
contained in the box $B_\infty(c^\epsilon_p)$, where
$$
c^\epsilon_p:=n\cdot\sup\limits_{1\leqslant i\leqslant
n}|(\Phi^{\top}(\Phi\Phi^{\top})^{-1}\mathbf{b}_{\epsilon})_{i}|.
$$
\end{Lemma}

\noindent {\bf Proof.}
This result holds from \cite[Remark 1]{peng15equivalence_Lp}, here
we omit the proof.
\ep

Now, we consider problem $(P_p)$ with $p\in (0,1)$.
\begin{Theorem}\label{PCS-Lp-implicitly-has-bounds}
Problem $(P_p)$  is equivalent to
\begin{eqnarray}\label{PCS_Pp-contain-bounds-implicitly}
\min \|\mathbf{x}\|_{p}^{p}\quad \textrm{s.t.}\;\;|\Phi\mathbf{x}|=\mathbf{b}, \|\mathbf{x}\|_{\infty}\leqslant c_1:=\max\limits_{\epsilon\in\{-1,1\}^m}c^\epsilon_p>0,
\end{eqnarray}
where $c^\epsilon_p$ is defined in Lemma \ref{CS-Lp-bounded}.
\end{Theorem}

\noindent {\bf Proof.}
Let $S_p$ be the solution set of problem $(P_p)$, then
$$
S_p\subseteq \bigcup\limits_{\epsilon\in\{-1,1\}^m} S^{p}_\epsilon.
$$
Since every $S^{p}_\epsilon$ is contained in the box $B_\infty(c^\epsilon_p)$ by
Lemma \ref{CS-Lp-bounded}, it follows from the above equality that $S_{p}$ is contained in $\bigcup_\epsilon B_\infty(c^\epsilon_p)\subseteq B_\infty(c_1)$ with $c_1$ being defined by (\ref{PCS_Pp-contain-bounds-implicitly}). Hence, the feasible set of problem (\ref{PCS_Pp-contain-bounds-implicitly}) contains all optimal solutions of problem $(P_p)$. By Lemma \ref{contains-all-solutions}, problem (\ref{PCS_Pp-contain-bounds-implicitly}) is equivalent to problem $(P_p)$. The proof is complete.
\ep

\subsection{Reformulations of Problems $(P_0)$ and $(P_p)$}

\hspace{4mm} Define a constant
\begin{eqnarray}\label{on-constant-c0}
r_0:=\max\{c_{0},c_{1}\}>0,
\end{eqnarray}
where $c_{0}, c_{1}$  are given in (\ref{PCS_P0-contain-bounds-implicitly}) and (\ref{PCS_Pp-contain-bounds-implicitly}), respectively.
\begin{Remark}\label{on-constant-c0-2}
From the results of the former two subsections, we can see that
\begin{eqnarray*}
(P_{0bd})\quad\min \|\mathbf{x}\|_{0}\quad \textrm{s.t.}\;\;|\Phi\mathbf{x}|=\mathbf{b}, \|\mathbf{x}\|_{\infty}\leqslant r_0
\end{eqnarray*}
and
\begin{eqnarray*}
(P_{pbd})\quad\min \|\mathbf{x}\|_{p}^{p}\quad \textrm{s.t.}\;\;|\Phi\mathbf{x}|=\mathbf{b}, \|\mathbf{x}\|_{\infty}\leqslant r_0
\end{eqnarray*}
are equivalent to problem $(P_0)$ and problem $(P_p)$, respectively.
\end{Remark}

Introduce a variable $\mathbf{y}\in\mathbb{R}^{n}$ such that $|\mathbf{x}|\leqslant\mathbf{y}$. Noting that $|\mathbf{x}|\leqslant \mathbf{y}$ implies that $-\mathbf{y}\leqslant\mathbf{x}\leqslant\mathbf{y}$, or equivalently,
$-\mathbf{x}-\mathbf{y}\leqslant\mathbf{0}$ and $\mathbf{x}-\mathbf{y}\leqslant\mathbf{0}$. Define two block matrices $\Psi:=\mathrm{\rm(-E ~,~E)^{\top}}$ and
$\Gamma:=\mathrm{\rm (-E , -E)^{\top}}$ where $\mathrm{\rm E}\in\mathbb{R}^{n\times n}$ is the $n\times n$ identity matrix, then it is easy to verify that
$|\mathbf{x}|\leqslant\mathbf{y}$ if and only if $\Psi \mathbf{x}+\Gamma \mathbf{y}\leqslant\mathbf{0}$.
Let
\begin{eqnarray*}
\mathbb{T}_{\epsilon}:=
\left\{(\mathbf{x}, \mathbf{y})\in \mathbb{R}^{2n}\left|
\begin{array}{ccc}
{\Phi \mathbf{x} =\mathbf{b}_{\epsilon},\Psi \mathbf{x}+\Gamma \mathbf{y}\leqslant\mathbf{0}}\\
{-r_0\mathbf{1}\leqslant \mathbf{x}\leqslant r_0\mathbf{1}}\\
{\mathbf{0}\leqslant \mathbf{y}\leqslant r_0\mathbf{1}}
\end{array}
\right.
\right\},
\end{eqnarray*}
where $\mathbf{\epsilon}\in\{-1,1\}^{m}\subset\mathbb{R}^{m}$ and $\mathbf{b}_{\epsilon}=\mathbf{b}\circ\epsilon$, then it holds that
\begin{eqnarray}\label{unionofpolytope}
\mathbb{T}:=\left\{(\mathbf{x},\mathbf{y})\in\mathbb{R}^{2n}\Bigg|
\begin{array}{ll}
{|\Phi\mathbf{x}|=\mathbf{b},\|\mathbf{x}\|_{\infty}\leqslant r_0}\\
{|\mathbf{x}|\leqslant\mathbf{y},\mathbf{0}\leqslant\mathbf{y}\leqslant r_0\mathbf{1}}
\end{array}
\right\}= \bigcup_{\mathbf{\epsilon}\in\{-1,1\}^{m}}\mathbb{T}_{\mathbf{\epsilon}}.
\end{eqnarray}
It is easy to see that
\begin{eqnarray*}
\mathbb{T}, \mathbb{T}_{\epsilon}\subseteq[-r_0,r_0]^n\times[0,r_0]^n\subset\mathbb{R}^{2n}.
\end{eqnarray*}

Now, we consider the following two problems with the same feasible set:
\begin{eqnarray}
  \min\limits_{(\mathbf{x},\mathbf{y})} \|\mathbf{y}\|_0\quad \textrm{s.t.}\;\; (\mathbf{x},\mathbf{y})\in \mathbb{T},\label{2.3} \\
 \min\limits_{(\mathbf{x},\mathbf{y})}\|\mathbf{y}\|_p^p\quad \textrm{s.t.}\;\; (\mathbf{x},\mathbf{y})\in \mathbb{T}.\label{lp-relax-original}
\end{eqnarray}

\begin{Theorem}\label{eqivalent-forms-of-p0-contain-bounds}
 Problem $ (P_{0bd})$ is equivalent to  problem {\rm (\ref{2.3})}, and problem $(P_{pbd})$ is equivalent to  problem {\rm (\ref{lp-relax-original})}.
\end{Theorem}

\noindent {\bf Proof.}
From Examples \ref{on-sign-fun} and \ref{on-lp-fun}, we know that both $\|\cdot\|_{0}$  and $\|\cdot\|_{p}^{p}$ are monotonically nondecreasing when they are confined from $\mathbb{R}^{n}_{+}$ to $\mathbb{R}_{+}$. Let
\begin{eqnarray*}
\begin{array}{l}
\hat{\mathbb{D}}_{R^{n}}:=\left\{\mathbf{x}\in\mathbb{R}^{n}\;\Big|\;|\Phi\mathbf{x}|=\mathbf{b}\right\},\\
\mathbb{D}_{R^{n}}:=\{\mathbf{x}\mid\mathbf{x}\in \hat{\mathbb{D}}_{R^n},\mathbf{l}\leqslant \mathbf{x}\leqslant \mathbf{u}\}, \\
\mathbb{D}_{R^{2n}}:=\{(\mathbf{x}, \mathbf{y})\mid\mathbf{x}\in \mathbb{D}_{R^{n}}, |\mathbf{x}|\leqslant \mathbf{y}, \mathbf{0}\leqslant \mathbf{y} \leqslant \mathbf{v}\},
\end{array}
\end{eqnarray*}
where $\mathbf{l}=- r_0\mathbf{1}<\mathbf{0}$ and $\mathbf{u}=r_0\mathbf{1}>\mathbf{0}$, and hence,
$\mathbf{v}=\max\{-\mathbf{l},\mathbf{u}\}=r_0\mathbf{1}$. By a direct application of Corollary \ref{add-cro1-1}, we complete the proof.
\ep

\begin{Remark}\label{all-equivalences}
From Remark \ref{on-constant-c0-2} and Theorem \ref{eqivalent-forms-of-p0-contain-bounds}, we have the following equivalent relationships:
$$
 {\rm (\ref{2.3})} \Longleftrightarrow(P_{0bd})\Longleftrightarrow(P_{0}),\qquad  {\rm (\ref{lp-relax-original})}\Longleftrightarrow(P_{pbd})\Longleftrightarrow(P_{p}).
$$

\noindent Hence, for some $p\in (0,1)$, if we can obtain the equivalence relationship between problem {\rm (\ref{2.3})} and problem {\rm (\ref{lp-relax-original})}, then we will achieve the equivalence relationship between problem $(P_{0})$ and problem $(P_{p})$.
\end{Remark}

\section{Equivalence of Problem ($P_0$) and Problem ($P_p$)}
\setcounter{equation}{0} \setcounter{Assumption}{0}
\setcounter{Theorem}{0} \setcounter{Proposition}{0}
\setcounter{Corollary}{0} \setcounter{Lemma}{0}
\setcounter{Definition}{0} \setcounter{Remark}{0}
\setcounter{Algorithm}{0} \setcounter{Example}{0}

\hspace{4mm} In this section, we study the equivalence relationship between problem $(P_p)$ and problem $(P_0)$. Some facts for the polytope are given in the following proposition which are useful for our subsequent discussions.
\begin{Proposition}\label{Representation theorem for polytopes}
\begin{itemize}
\item[(i)] \textrm{\rm (\cite[Theorem 2.15(6)]{Ziegler95bk})}
A subset $\mathbb{P}\subseteq \mathbb{R}^d$ is a polytope if and only if it can be described as a bounded intersection of closed half-spaces.
\item[(ii)] \textrm{\rm (\cite[Theorem 1.1]{Ziegler95bk})}
A subset $\mathbb{P}\subseteq\mathbb{R}^d$ is the convex hull of a finite point set: $\mathbb{P}=\textrm{Conv}(\mathbb{V})$ for some $\mathbb{V}\in\mathbb{R}^{d\times n}$
if and only if it is a bounded intersection of half-spaces.
\item[(iii)] \textrm{\rm (\cite[Proposition 2.2]{Ziegler95bk})}
Let $\mathbb{P}\subseteq \mathbb{R}^d$ be a polytope.
(a) Every polytope is the convex hull of its vertices: $\mathbb{P}=\mathrm{conv(vert}{(\mathbb{P})})$.
(b) If a polytope can be written as the convex hull of a finite point set, then the set contains all the vertices of the polytope:
$\mathbb{P}=\mathrm{conv}(\mathbb{V})$ implies that $\mathrm{vert(\mathbb{P})}\subseteq \mathbb{V}$.
\item[(iv)] \textrm{\rm (\cite[ page32, item 4]{Grunbaum03bk})}
(a) The convex hull of finitely many polytopes is a polytope.
(b) The intersection of a polytope with an affine variety is a polytope.
\end{itemize}
\end{Proposition}

We now apply theory of polytope to investigate the set $\mathbb{T}$ defined by (\ref{unionofpolytope}).

\begin{Lemma}\label{about_T}
The set $\mathbb{T}$ defined by (\ref{unionofpolytope}) is non-empty  and it is  the union  of finitely many
polytopes in $\mathbb{R}^{2n}$.
\end{Lemma}

\noindent {\bf Proof.}
Since we have assumed in problem $(P_0)$ that $m\ll n$ and $\Phi$ has full row rank, it follows that the feasible set of problem $(P_{0})$ is non-empty. Combining
the equivalence relationships in Remark \ref{all-equivalences}, we see
that $~\mathbb{T}$ is a non-empty set.  By \reff{unionofpolytope},
we only need to prove $~\mathbb{T}_{\epsilon}~$ is a polytope. This
can be seen from that it is a bounded intersection of finitely
closed half spaces and  an affine set (i.e., point set satisfying
equality constraints) in $\mathbb{R}^{2n}$. By (i) and (iv)(b) of Proposition
\ref{Representation theorem for polytopes}, $\mathbb{T}_{\epsilon}$
is a polytope in $\mathbb{R}^{2n}$.
\ep

It should be noticed that the set $\mathbb{T}$  is  not  convex, since it
is a union set of disjoint convex sets (i.e, polytopes).

By Lemma \ref{about_T}, problem ($P_0$) is equivalent to a
sparse optimization problem over a  union set of finitely many disjoint polytopes.
As is well known that the objective function of problem \reff{2.3} is
combinatorial and the feasible set is non-convex, so problem
\reff{2.3} is NP-hard \cite{Nata95linearL0-NP-hard}.

Define the convex hull of the set $\mathbb{T} $ given in (\ref{unionofpolytope}) by
\begin{eqnarray}\label{2.4}
\mathbb{S}:= \textrm{Conv}(\mathbb{T}) =
\textrm{Conv}\left(\bigcup_{\epsilon\in\{-1,1\}^{m}}\mathbb{T}_{\epsilon}\right)\subseteq [-r_0,r_0]^n\times[0,r_0]^n\subset\mathbb{R}^{2n}.
\end{eqnarray}
Then, we have the following result.
\begin{Lemma}\label{OnS}
The set $\mathbb{S}$ defined by (\ref{2.4}) is a non-empty polytope contained in the box $[-r_0,r_0]^n\times[0,r_0]^n\subset\mathbb{R}^{2n}$.
\end{Lemma}

\noindent {\bf Proof.} By Lemma \ref{about_T}, we know that the set $\mathbb{T}$ is non-empty and  is  the union  of finitely many
polytopes in $\mathbb{R}^{2n}$. Thus, the set $\mathbb{S}$ is  non-empty and is a convex hull of finitely many polytopes. By
(iv)(a) of Proposition \ref{Representation theorem for polytopes}, the set
$\mathbb{S}$ is a polytope. Furthermore, $\mathbb{S}$ is the
smallest convex set containing
$\mathbb{T}\subset[-r_0,r_0]^n\times[0,r_0]^n$.
\ep

For the set $\mathbb{T}$ defined by (\ref{unionofpolytope}), we call $\mathbf{x}\in \mathbb{T}$ a pseudo-extreme point of $\mathbb{T}$ if it is a vertex of some polytope $\mathbb{T}_{\epsilon}$, which implies that there are no two distinct points $\mathbf{y}, \mathbf{z}\in \mathbb{T}_{\epsilon}$ and a constant $\lambda\in(0,1)$ such that $\mathbf{x}=\lambda \mathbf{y}+(1-\lambda)\mathbf{z}$;
an extreme point of $\mathbb{S}$  is called a vertex of $\mathbb{S}$ as usually defined. In the following, we use $\mathrm{\rm vert }(\mathbb{S})$ to denote the set of all vertices of $\mathbb{S}$ and $\mathrm{\rm pext}(\mathbb{T})$ to denote the set of all pseudo-extreme points of $\mathbb{T}$.

Using the convex hull $\mathbb{S}$, problem (\ref{lp-relax-original}) is relaxed as
\begin{eqnarray}\label{2.5}
\min\limits_{(\mathbf{x},\mathbf{y})}\|\mathbf{y}\|_p^p\quad \textrm{s.t.}\; (\mathbf{x},\mathbf{y})\in \mathbb{S}.
\end{eqnarray}
We have the following result.
\begin{Theorem}\label{add-lem2-2}
Suppose that $\mathbf{x}^*\in \mathrm{\rm vert }(\mathbb{S})$ is an optimal solution to  problem {\rm(\ref{2.5})}, which is called a vertex solution of problem
{\rm(\ref{2.5})}, then $\mathbf{x}^*\in \mathrm{\rm pext}(\mathbb{T})$ and it solves problem {\rm(\ref{lp-relax-original})}, which is called a pseudo-extreme
point solution of problem {\rm(\ref{lp-relax-original})}.
\end{Theorem}

\noindent {\bf Proof.}
Notice that problems (\ref{2.5}) and (\ref{lp-relax-original}) share the same objective function while have two different feasible sets, i.e.,  $\mathbb{T}$ and $\mathbb{S}$ with $\mathbb{T}\subset\mathbb{S}$. Since $\mathbb{S}=\textrm{Conv}(\mathrm{\rm pext}(\mathbb{T}))$, it follows from (iii)(b) of Proposition \ref{Representation theorem for polytopes} that $\mathrm{\rm vert }(\mathbb{S})\subseteq\mathrm{\rm pext}(\mathbb{T}$). Therefore, by Lemma \ref{add-lem2-1}, each vertex solution to problem (\ref{2.5}) is a pseudo-extreme point solution to problem (\ref{lp-relax-original}).
\ep

The following proposition will be used.
\begin{Proposition}\label{Rock70bk_Corollary2_P345}\textrm{\rm (\cite[Corollary 32.3.4]{Rock70bk})}
Let $f$ be a convex function, and let $C$ be a non-empty polyhedral convex set contained in $dom~f$. Suppose that  $C$ contains no lines,  and that $f$ is bounded above on $C$. Then the supremum of $f$
relative to $C$ is attained   at one of the (finitely many) extreme points of $C$.
\end{Proposition}

\begin{Theorem}\label{add-lem3-1} For each real number $p\in (0,1)$, there exists a pseudo-extreme point solution to problem (\ref{lp-relax-original}).
 \end{Theorem}

\noindent {\bf Proof.}
Denote $f(\mathbf{x},\mathbf{y}):=\|\mathbf{y}\|_p^p=\sum\limits_{i=1}^ny_i^p$ with $p\in (0,1)$, then $-f$ is a convex function on
$\textrm{dom}f=\mathbb{R}^n\times\mathbb{R}^n_+$. By Lemma \ref{OnS}, the set $\mathbb{S}$ is non-empty and
$\mathbb{S}\subseteq[-r_0,r_0]^n\times[0,r_0]^n\subset\textrm{dom}f$. There are no half-lines and no lines in $\mathbb{S}$; and
$-f$ is bounded above on $\mathbb{S}$. By Proposition \ref{Rock70bk_Corollary2_P345}, there exists a vertex solution to problem {\rm(\ref{2.5})}, which is also a pseudo-extreme point solution to problem (\ref{lp-relax-original}) by Theorem \ref{add-lem2-2}. In other words, for each real number $p\in (0, 1)$, there exists a pseudo-extreme point solution to problem (\ref{lp-relax-original}).
\ep

\begin{Theorem}\label{Pp-solution-is-extreme-point-solution-T}
Suppose that $\mathbf{x}_p$ is an optimal solution of problem {\rm($P_p$)}. Then, $(\mathbf{x}_p,|\mathbf{x}_p|)$ is a pseudo-extreme point solution of problem {\rm (\ref{lp-relax-original})}.
\end{Theorem}

\noindent {\bf Proof.}
Since $\mathbf{x}_p$ is an optimal solution of problem {\rm($P_p$)}, it follows from the equivalence relationship $(P_{pbd})\Leftrightarrow(P_{p})$ given in Remark \ref{all-equivalences} that $\mathbf{x}_p\in B_\infty(r_0)$ and $\mathbf{x}_p$ is also an optimal solution of problem {\rm($P_{pbd}$)}. Moreover, we can obtain that $(\mathbf{x}_p,\mathbf{y}_p)$ with $\mathbf{y}_p=|\mathbf{x}_p|$ is an optimal solution of problem {\rm (\ref{lp-relax-original})}. In fact, for any $(\mathbf{x},\mathbf{y})\in\mathbb{T}$, by noticing that $|\mathbf{x}|\leqslant\mathbf{y}$ and the nondecreasing property of $\|\cdot\|_{p}^{p}$, we have
$$
\|\mathbf{y}_p\|_{p}^{p}=\||\mathbf{x}_p|\|_{p}^{p} =\|\mathbf{x}_p\|_{p}^{p}\leqslant\|\mathbf{x}\|^{p}_{p}\leqslant\|\mathbf{y}\|^{p}_{p},
$$
which implies that $(\mathbf{x}_p,|\mathbf{x}_p|)$ solves problem {\rm (\ref{lp-relax-original})}.

Next, we show that $(\mathbf{x}_p,|\mathbf{x}_p|)$ is a pseudo-extreme point of $\mathbb{T}$ by contradiction.
Since a pseudo-extreme point of $\mathbb{T}$ is a vertex  of polytope $\mathbb{T}_\epsilon$ for some $\epsilon\in\{-1,1\}^m$, we assume that $(\mathbf{x}_p,|\mathbf{x}_p|)\in\mathbb{T}_\epsilon$ with any fixed $\epsilon\in\{-1,1\}^m$ is not a vertex of $\mathbb{T}_\epsilon$, then there exist some fixed $\lambda\in(0,1)$ and two distinct points $(\mathbf{x}^1, \mathbf{y}^1)$ and $(\mathbf{x}^2, \mathbf{y}^2)$ in $\mathbb{T}_\epsilon$ with $\mathbf{x}^1\neq\mathbf{x}^2$ or $\mathbf{y}^1\neq\mathbf{y}^2$ such that
\begin{eqnarray}\label{convex-combination}
(\mathbf{x}_p,|\mathbf{x}_p|)
=\lambda(\mathbf{x}^1, \mathbf{y}^1)+(1-\lambda)(\mathbf{x}^2, \mathbf{y}^2)
=\left(\lambda\mathbf{x}^1+(1-\lambda)\mathbf{x}^2,\lambda\mathbf{y}^1+(1-\lambda)\mathbf{y}^2\right).
\end{eqnarray}
By $(\mathbf{x}^1, \mathbf{y}^1),(\mathbf{x}^2, \mathbf{y}^2)\in \mathbb{T}_\epsilon$, we have
\begin{eqnarray}\label{inter-result-1}
\mathbf{0}\leqslant|\mathbf{x}^1|\leqslant\mathbf{y}^1\leqslant r_0\mathbf{1}\quad \mbox{\rm and}\quad \mathbf{0}\leqslant|\mathbf{x}^2|\leqslant\mathbf{y}^2\leqslant r_0\mathbf{1}.
\end{eqnarray}

Case (i). If $\mathbf{y}^1=\mathbf{y}^2$, we have  $\mathbf{x}^1\neq\mathbf{x}^2$ and $|\mathbf{x}_p|=\mathbf{y}^1=\mathbf{y}^2$ from (\ref{convex-combination}). By (\ref{convex-combination}) and (\ref{inter-result-1}), we have
\begin{eqnarray}\label{inter-result-2}
|\mathbf{x}_p|=|\lambda\mathbf{x}^1+(1-\lambda)\mathbf{x}^2|\leqslant\lambda|\mathbf{x}^1|+(1-\lambda)|\mathbf{x}^2|
\leqslant\lambda\mathbf{y}^1+(1-\lambda)\mathbf{y}^2=|\mathbf{x}_p|.
\end{eqnarray}
From (\ref{inter-result-2}), it follows that $$|\lambda\mathbf{x}^1+(1-\lambda)\mathbf{x}^2|=\lambda|\mathbf{x}^1|+(1-\lambda)|\mathbf{x}^2|,$$ which implies $\mathbf{x}^1_i\mathbf{x}^2_i\geqslant 0$ for all $i\in [n]$. Therefore,
\begin{eqnarray}\label{inter-result-3}
\|\mathbf{x}_p\|_{p}^{p}&=&\left\||\mathbf{x}_p|\right\|_{p}^{p}=\left\||\lambda\mathbf{x}^1+(1-\lambda)\mathbf{x}^2|\right\|_{p}^{p}\nonumber\\ &=&\left\|\lambda|\mathbf{x}^1|+(1-\lambda)|\mathbf{x}^2|\right\|_{p}^{p}\nonumber\\
&>&\lambda\left\||\mathbf{x}^1|\right\|_{p}^{p}+(1-\lambda)\left\||\mathbf{x}^2|\right\|_{p}^{p}\nonumber\\
&=&\lambda\|\mathbf{x}^1\|_{p}^{p}+(1-\lambda)\|\mathbf{x}^2\|_{p}^{p}
\end{eqnarray}
where the strict inequality is due to the strict concavity of $\|\cdot\|_{p}^{p}$ on
$\mathbb{R}^n_{+}$ when $p\in (0,1)$. Without loss of generality, let $\|\mathbf{x}^1\|_{p}^{p}=\min\{\|\mathbf{x}^1\|_{p}^{p},\|\mathbf{x}^2\|_{p}^{p}\}$. By {\rm(\ref{inter-result-3})}, $\|\mathbf{x}_p\|_{p}^{p}>\|\mathbf{x}^1\|_{p}^{p}$. This contradicts the condition that $\mathbf{x}_p$
is an optimal solution of {\rm($P_p$)}.

Case (ii). We assume $\mathbf{y}^1\neq\mathbf{y}^2$. By (\ref{convex-combination}),
we have $|\mathbf{x}_p|=\lambda\mathbf{y}^1+(1-\lambda)\mathbf{y}^2$. Furthermore,
\begin{eqnarray}\label{inter-result-4}
\|\mathbf{x}_p\|_{p}^{p}&=&\left\||\mathbf{x}_p|\right\|_{p}^{p}=\|\lambda\mathbf{y}^1+(1-\lambda)\mathbf{y}^2\|^{p}_{p}\nonumber\\
                        &>&\lambda\|\mathbf{y}^1\|^{p}_{p}+(1-\lambda)\|\mathbf{y}^2\|^{p}_{p}\nonumber\\
                        &\geqslant&\lambda\|\mathbf{x}^1\|^{p}_{p}+(1-\lambda)\|\mathbf{x}^2\|^{p}_{p},
\end{eqnarray}
where the first inequality is due to  the strict concavity of $\|\cdot\|_{p}^{p}$ on
$\mathbb{R}^n_{+}$ when $p\in (0,1)$, and the second inequality is due to $\|\cdot\|_p^p$ being nondecreasing on $\mathbb{R}_+^n$. Without loss of generality, let $\|\mathbf{x}^1\|_{p}^{p}=\min\{\|\mathbf{x}^1\|_{p}^{p},\|\mathbf{x}^2\|_{p}^{p}\}$. By {\rm(\ref{inter-result-4})}, $\|\mathbf{x}_p\|_{p}^{p}>\|\mathbf{x}^1\|_{p}^{p}$. This contradicts the condition that $\mathbf{x}_p$
is an optimal solution of problem {\rm($P_p$)}.
From the above discussions, we see that $(\mathbf{x}_p,|\mathbf{x}_p|)$ can not be expressed as a convex combination of two distinct points of any fixed $\mathbb{T}_\epsilon$.
Hence, $(\mathbf{x}_p,|\mathbf{x}_p|)$ is a pseudo-extreme point of $\mathbb{T}$.

The proof is complete.
\ep

Define
\begin{eqnarray}\label{min-positive-coordinate-of-extreme-point}
r_m:=\min\limits_{(\tilde{\mathbf{x}},\tilde{\mathbf{y}})\in{\rm pext\left(\mathbb{T}\right)},\tilde{y}_k\neq 0 }{\tilde{y}_k}.
\end{eqnarray}
By (ii) and (iii)(a) of Proposition \ref{Representation theorem for polytopes}, the number of vertices of any polytope is finite. Thus, as the union of finitely many disjoint polytopes, $\mathbb{T}$ has finitely many pseudo-extreme points. Hence we have $ r_m>0$. Moreover, by noticing that $\mathbf{0}\leqslant\tilde{\mathbf{y}}\leqslant r_0\mathbf{1}$, where $r_0$ is defined by (\ref{on-constant-c0}), we can obtain that $r_0\geqslant r_m$.

Recall the equivalence relationships in Remark \ref{all-equivalences}, that is,
{\rm (\ref{2.3})} $\Leftrightarrow(P_{0bd})\Leftrightarrow(P_{0})$ and $(P_{p})\Leftrightarrow(P_{pbd})\Leftrightarrow$ {\rm (\ref{lp-relax-original})},
we have
\begin{eqnarray}\label{minL0-bounded-equals-minL0}
\min\limits_{|\Phi\mathbf{x}|=\mathbf{b},\mathbf{x}\in B_{\infty}(r_0)}\|\mathbf{x}\|_{0}=\min\limits_{|\Phi\mathbf{x}|=\mathbf{b}}\|\mathbf{x}\|_{0},
\end{eqnarray}
\begin{eqnarray}\label{minLp-bounded-equals-minLp}
 \min\limits_{{(\mathbf{x},\mathbf{y})\in{\rm pext\left(\mathbb{T}\right)}}}\|\mathbf{y}\|^{p}_{p}=\min\limits_{(\mathbf{x},\mathbf{y})\in{\rm \mathbb{T}}}\|\mathbf{y}\|^{p}_{p}=\min\limits_{|\Phi\mathbf{x}|=\mathbf{b},\mathbf{x}\in B_{\infty}(r_0)}\|\mathbf{x}\|^{p}_{p},
\end{eqnarray}
where the first equality of (\ref{minLp-bounded-equals-minLp}) is due to Theorem \ref{add-lem3-1}.

Now, in a similar way as in \cite{peng15equivalence_Lp}, we construct a constant $p^\ast\in (0,1]$ such that for any fixed $p\in (0,p^\ast)$, an arbitrarily optimal solution to problem $(P_p)$, say ${\mathbf{x}}_p$, is also an optimal solution to problem $(P_0)$.

Since ${\mathbf{x}}_p$ is an optimal solution to problem $(P_p)$, by virtue of Theorem {\rm\ref{Pp-solution-is-extreme-point-solution-T}}, it follows that $({\mathbf{x}}_p,{\mathbf{y}}_p)$ with ${\mathbf{y}}_p=|{\mathbf{x}}_p|$ is a pseudo-extreme point solution to problem (\ref{lp-relax-original}).
From the definition of $r_m$ in (\ref{min-positive-coordinate-of-extreme-point}), we have $\left(\frac{{(y_p)}_i}{r_m}\right)^{p}\geqslant 1$.
Furthermore,
\begin{eqnarray}\label{find-pstar}
 \|{\mathbf{x}}_p\|_{0}
 &=&\|{\mathbf{y}}_p\|_{0} =\lim\limits_{p\downarrow0}\sum\limits_{i=1}^{n}\left(\frac{{(y_p)}_i}{r_m}\right)^{p}
 \leqslant\sum\limits_{i=1}^{n}\left(\frac{{(y_p)}_i}{r_m}\right)^{p}\nonumber\\
&=& r_m^{-p}\min\limits_{{(\mathbf{x},\mathbf{y})\in{\rm pext\left(\mathbb{T}\right)}}}\|\mathbf{y}\|^{p}_{p}
= r_m^{-p}\min\limits_{|\Phi\mathbf{x}|=\mathbf{b},\mathbf{x}\in B_{\infty}(r_0)}\|\mathbf{x}\|^{p}_{p}\nonumber\\
&=&\left(\frac{r_0}{r_m}\right)^{p} \min\limits_{|\Phi\mathbf{x}|=\mathbf{b},\mathbf{x}\in B_{\infty}(r_0)}\|r_0^{-1}\mathbf{x}\|^{p}_{p}\nonumber\\
&\leqslant &\left (\frac{r_0}{r_m}\right )^{p} \min\limits_{|\Phi\mathbf{x}|=\mathbf{b},\mathbf{x}\in B_{\infty}(r_0)}\|r_0^{-1}\mathbf{x}\|_{0}\nonumber\\
&=&\left(\frac{r_0}{r_m}\right)^{p}\min\limits_{|\Phi\mathbf{x}|=\mathbf{b},\mathbf{x}\in B_{\infty}(r_0)}\|\mathbf{x}\|_{0}\nonumber\\
&=&\left(\frac{r_0}{r_m}\right)^{p}\min\limits_{|\Phi\mathbf{x}|=\mathbf{b}}\|\mathbf{x}\|_{0},
\end{eqnarray}
where the first inequality follows from the nondecreasing property of function $f(p)=t^p$ when $p\in(0,1)$ with each fixed $ t\geqslant1$,
the second inequality from the non-increasing property of function $f(p)=t^p$ when $p\in [0,1)$ with each fixed $0\leqslant t<1$ and $0^p\equiv0$,
the fourth equality from (\ref{minLp-bounded-equals-minLp}), and
the last equality from (\ref{minL0-bounded-equals-minL0}).

Since $\|{\mathbf{x}}_p\|_{0}$ is an integer, it follows from (\ref{find-pstar}) that  ${\mathbf{x}}_p$ is an optimal solution of problem $(P_0)$ when
\begin{eqnarray}\label{use-this-expression-to-decide-pStar-1}
\|{\mathbf{x}}_p\|_{0}\leqslant \left (\frac{r_0}{r_m}\right)^{p}\min\limits_{|\Phi\mathbf{x}| =\mathbf{b}}\|\mathbf{x}\|_{0}<\min\limits_{|\Phi\mathbf{x}|=\mathbf{b}}\|\mathbf{x}\|_{0}+1.
\end{eqnarray}
We assume that $\min\limits_{|\Phi\mathbf{x}|=\mathbf{b}}\|\mathbf{x}\|_{0}=s\geqslant1$. Then, it follows from (\ref{use-this-expression-to-decide-pStar-1}) that ${\mathbf{x}}_p$ solves problem $(P_0)$ when $p\in (0,1)$ satisfies
\begin{eqnarray}\label{use-this-expression-to-decide-pStar-2}
s\leqslant  \left(\frac{r_0}{r_m}\right)^{p}s<s+1.
\end{eqnarray}
Two cases are discussed as follows.
\begin{itemize}
  \item [(i)] $\mbox{\rm If}\;{1\leqslant\frac{r_0}{r_m}\leqslant\frac{s+1}{s}}$, then  ${s\leqslant\frac{r_0}{r_m}s\leqslant s+1}$. We conclude  ${s\leqslant\left(\frac{r_0}{r_m}\right)^{p}s<s+1}$ holds for all $p\in (0,1)$. In this case, we can define $p^\ast:=1$.
  \item [(ii)] $\mbox{\rm If}\;{1<\frac{s+1}{s}<\frac{r_0}{r_m}}$, then  ${s+1<\frac{r_0}{r_m}s}$. In this case, we need adapt $p\in (0,1)$  to satisfy (\ref{use-this-expression-to-decide-pStar-2}). By taking logarithm for each term of (\ref{use-this-expression-to-decide-pStar-2}), we derive $\ln s \leqslant p(\ln r_0 -\ln r_m)+\ln s<\ln (s+1)$,
      and hence, $0 \leqslant p(\ln r_0 -\ln r_m)<\ln (s+1)-\ln s$. Therefore
     $$0 \leqslant p<\frac{\ln (s+1)-\ln s}{\ln r_0 -\ln r_m}.$$
     So, in this case, we can define $p^\ast:=\frac{\ln (s+1)-\ln s}{\ln r_0 -\ln r_m}$.
\end{itemize}

With the help of the above discussions, we obtain the following result immediately.
\begin{Theorem}\label{mainresult-2}
Suppose that $\min\limits_{|\Phi\mathbf{x}|=\mathbf{b}}\|\mathbf{x}\|_{0}=s\geqslant1$.
Let
\begin{eqnarray}\label{p-star-formula}
p^\ast :=
\left\{
\begin{array}{cl}
{1}, & \mbox{\rm if}\;{1\leqslant\frac{r_0}{r_m}\leqslant\frac{s+1}{s}},\\
{\frac{\ln(s+1)-\ln s}{\ln r_0-\ln r_m}}, & \mbox{\rm if}\;{\frac{s+1}{s}<\frac{r_0}{r_m}},
\end{array}
\right.
\end{eqnarray}
where $r_0$ and $r_m$ are defined by (\ref{on-constant-c0}) and (\ref{min-positive-coordinate-of-extreme-point}), respectively.
Then,  for any fixed $p\in (0,p^\ast)$, every optimal solution to problem $(P_p)$ solves problem $(P_0)$.
\end{Theorem}

\begin{Remark}
Obviously, $p^\ast$ is not unique, since any fixed $p$ with $0<p<p^\ast$ can play the same role as the earlier found $p^\ast$. Moreover, the value of $p^\ast$ constructed in (\ref{p-star-formula}) depends on the sparsity level of the concerned problem.
\end{Remark}

\section{Examples}
\setcounter{equation}{0} \setcounter{Assumption}{0}
\setcounter{Theorem}{0} \setcounter{Proposition}{0}
\setcounter{Corollary}{0} \setcounter{Lemma}{0}
\setcounter{Definition}{0} \setcounter{Remark}{0}
\setcounter{Algorithm}{0} \setcounter{Example}{0}

\hspace{4mm} In this section, we give two simple examples to illustrate our theoretical findings.

\begin{Example}\label{example1-lp-solve-sparse}
Consider a 2-dimensional problem:
\begin{eqnarray}\label{e-exam1-1}
\min\limits_{x\in \mathbb{R}^2}\| x \|_0\quad \textrm{s.t.}\;\; |x_1-x_2|=1.
\end{eqnarray}
\end{Example}

\begin{figure}[ht]
\centering
\includegraphics[width=2.5in]{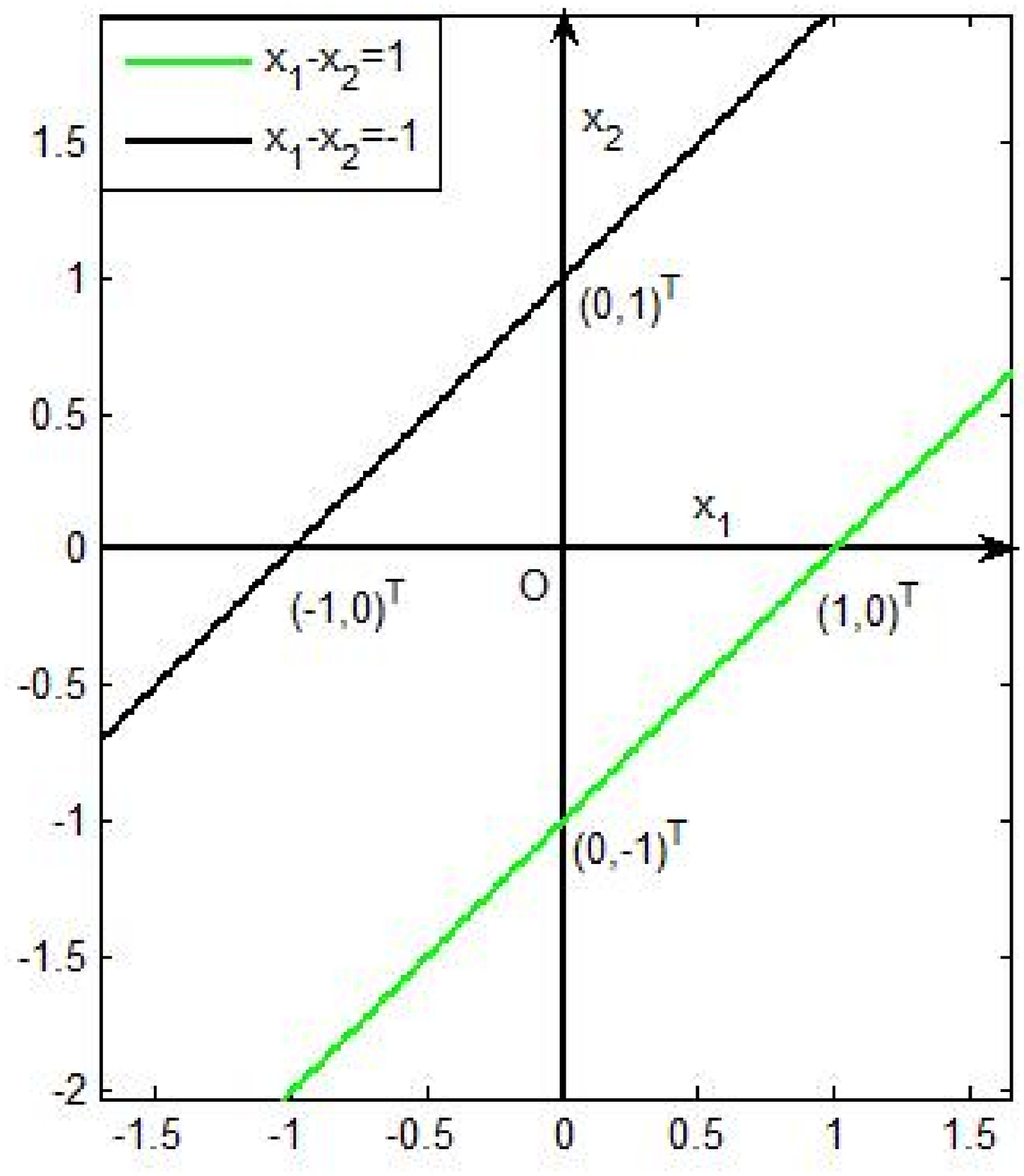}
\caption{The problem in Example \ref{example1-lp-solve-sparse} has four sparse solutions.}
\label{examples1-1}
\end{figure}

First, it is easy to verify that this problem has four sparse solutions: $\hat{\mathbf{x}}^1:=(0,-1)^\top$, $\hat{\mathbf{x}}^2:=(0,1)^\top$, $\hat{\mathbf{x}}^3:=(1,0)^\top$ and $\hat{\mathbf{x}}^4:=(-1,0)^\top$ (See Figure \ref{examples1-1}).

Second, we compute the value of $p^\ast$. It is easy to see that the sparsity level of the $\ell_0$-minimization (\ref{e-exam1-1}) is $s=1$; $\Phi=(1,-1)\in\mathbb{R}^{1\times2}$; $\epsilon\in\{-1,1\}$ and $\mathbf{b}=1$.
Denote $\epsilon_1:=-1$ and $\epsilon_2:=1$, then  $\mathbf{b}_{\epsilon_{1}}=\mathbf{b}\circ\epsilon_{1}=-1$ and $\mathbf{b}_{\epsilon_{2}}=\mathbf{b}\circ\epsilon_{2}=1$. Since
\begin{eqnarray*}
\begin{array}{rcl}
c^{\epsilon_1}_0&=&\max\limits_{I\subset[2],\#(I)=1}\max\limits_{i\in I}|((\Phi^{\top}_{I}\Phi_{I})^{-1}\Phi^{\top}_{I}\mathbf{b}_{\epsilon_1})_i| = \max\{1,1\}=1,\\  c^{\epsilon_2}_0&=&\max\limits_{I\subset[2],\#(I)=1}\max\limits_{i\in I}|((\Phi^{\top}_{I}\Phi_{I})^{-1}\Phi^{\top}_{I}\mathbf{b}_{\epsilon_2})_i|=\max\{1,1\}=1,
\end{array}
\end{eqnarray*}
it follows that
$$
c_0=\max\limits_{\epsilon\in\{-1,1\}}c^\epsilon_0=\max\{c^{\epsilon_1}_0,c^{\epsilon_2}_0\}=1;
$$
and since
\begin{eqnarray*}
\begin{array}{rcl}
c^{\epsilon_1}_p&=&n\cdot\sup\limits_{1\leqslant i\leqslant
n}|(\Phi^{\top}(\Phi\Phi^{\top})^{-1}\mathbf{b}_{\epsilon_1})_{i}| \\ &=&2\cdot\sup\limits_{i=1,2}|\left((1,-1)^{\top}[(1,-1)(1,-1)^{\top}]^{-1}\cdot(-1)\right)_i| =1,\\
c^{\epsilon_2}_p&=&n\cdot\sup\limits_{1\leqslant i\leqslant
n}|(\Phi^{\top}(\Phi\Phi^{\top})^{-1}\mathbf{b}_{\epsilon_2})_{i}| \\ &=&2\cdot\sup\limits_{i=1,2}|\left((1,-1)^{\top}[(1,-1)(1,-1)^{\top}]^{-1}\cdot1\right)_i| =1,
\end{array}
\end{eqnarray*}
it follows that
$$
c_1=\max\{c^{\epsilon_1}_p,c^{\epsilon_2}_p\}=\max\{1,1\}=1.
$$
Thus, $r_0=\max\{c_0,c_1\}=\max\{1,1\}=1$.
It is easy to see that the set $\mathbb{T}$ has four pseudo-extreme points $(0, -1, 0, 1)^{\top}$, $(0, 1, 0, 1)^{\top}$, $(1, 0, 1, 0)^{\top}$, $(-1, 0, 1, 0)^{\top}$, and hence, $r_m=1$. Furthermore,
$$2=\frac{s+1}{s}>\frac{r_0}{r_m}=\frac{1}{1}=1.$$
So, it follows from (\ref{p-star-formula}) that $p^\ast=1$.

Third, for any fixed $p\in (0,p^\ast)$, we consider the corresponding $\ell_p$-minimization:
\begin{eqnarray}\label{e-exam1-2}
\min\limits_{x\in \mathbb{R}^2}\| x \|_p^p\quad \textrm{s.t.}\;\; |x_1-x_2|=1,
\end{eqnarray}
and find all optimal solutions of this problem.

The $\ell_p$-minimization (\ref{e-exam1-2}) is equivalent to
$$
\min\{\underbrace{\min\limits_{x_1-x_2=1}\| \mathbf{x} \|_p^p}\limits_{(P_{p_{sub1}})},\;\underbrace{\min\limits_{x_1-x_2=-1}\| \mathbf{x } \|_p^p}\limits_{(P_{p_{sub2}})}\}.
$$
Both sub-problems can be equivalently transformed into one-variable unconstrained optimization problems, i.e.,
\begin{eqnarray}\label{exam1-subp}
(\tilde{P}_{p_{sub1}}):\;\;\min\limits_{x_1\in\mathbb{R}}|x_1|^p+|x_1-1|^p\quad
\mbox{\rm and}\quad
(\tilde{P}_{p_{sub2}}):\;\;\min\limits_{x_1\in\mathbb{R}}|x_1|^p+|x_1+1|^p.
\end{eqnarray}

\begin{figure}[ht]
\centering
\includegraphics[width=2.5in]{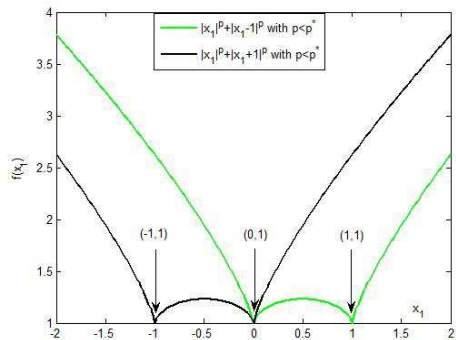}
\caption{Graphs of objective functions of two subproblems in (\ref{exam1-subp}) with one-variable $x_1$.}
\label{examples1-2}
\end{figure}

For sub-problem $(\tilde{P}_{p_{sub1}})$, we divide $\mathbb{R}$ into three parts: $\mathbb{R}=\mathbb{D}_1\bigcup\mathbb{D}_2\bigcup\mathbb{D}_3$
where $\mathbb{D}_1=(-\infty,0]$, $\mathbb{D}_2=[0,1]$ and $\mathbb{D}_3=[1,+\infty)$; and accordingly, we denote $f_1(x_1)=(-x_1)^p+(1-x_1)^p,f_2(x_1)=x_1^p+(1-x_1)^p$ and $f_3(x_1)=x_1^p+(x_1-1)^p$. Then, problem $(\tilde{P}_{p_{sub1}})$ is equivalent to
$$
\min\{\underbrace{\min\limits_{x_1\in\mathbb{D}_1}f_1(x_1)}\limits_{(\tilde{P}_{p^1_{sub1}})},
\underbrace{\min\limits_{x_1\in\mathbb{D}_2}f_2(x_1)}\limits_{(\tilde{P}_{p^2_{sub1}})},
\underbrace{\min\limits_{x_1\in\mathbb{D}_3}f_3(x_1)}\limits_{(\tilde{P}_{p^3_{sub1}})}\}.
$$
In the following, we consider the above three sub-problems, respectively.
\begin{itemize}
  \item [(i)] For problem $(\tilde{P}_{p^1_{sub1}})$. Since $f^\prime_1(x_1)=(-1)p(-x_1)^{p-1}-p(1-x_1)^{p-1}<0$ on $\mathbb{D}_1$, it follows that the function $f_1(\cdot)$ is strictly decreasing on $\mathbb{D}_1$. Thus, we have
      $\min\limits_{x_1\in\mathbb{D}_1}f_1(x_1)=f_1(0)=1^p=1$.
 \item [(ii)] For problem $(\tilde{P}_{p^2_{sub1}})$. It is easy to see that on $\mathbb{D}_2$,
     \begin{eqnarray*}
     \begin{array}{l}
     f^\prime_2(x_1)=px_1^{p-1}-p(1-x_1)^{p-1},\\ f^{\prime\prime}_2(x_1)=p(p-1)x_1^{p-2}+p(p-1)(1-x_1)^{p-2}<0.
     \end{array}
     \end{eqnarray*}
     Thus, the function $f_2(\cdot)$ is strictly concave on $\mathbb{D}_2$, which implies that  the optimal value is obtained on some end-point of $\mathbb{D}_2$. So,
     $$
     \min\limits_{x_1\in\mathbb{D}_2}f_2(x_1)=\min\{f_2(0),f_2(1)\}=f_2(0)=1^p=1.
     $$
\item [(iii)] For problem $(\tilde{P}_{p^3_{sub1}})$. Since $f'_3(x_1)=px_1^{p-1}+p(x_1-1)^{p-1}>0$
    on $\mathbb{D}_3$, it follows that the function $f_3(\cdot)$ is strictly increasing on $\mathbb{D}_3$. Thus,
      $\min\limits_{x_1\in\mathbb{D}_3}f_3(x_1)=f_3(1)=1^p=1$.
\end{itemize}
Combining cases (i)-(iii), we obtain that $\hat{x}^1_1=0$ and $\hat{x}^3_1=1$ are two optimal solutions to problem $(\tilde{P}_{p_{sub1}})$ (see the green curve in Figure \ref{examples1-2}). Furthermore, $\hat{\mathbf{x}}^1=(\hat{x}^1_1,\hat{x}^1_2)^\top=(0,-1)^\top$ and $\hat{\mathbf{x}}^3=(\hat{x}^3_1,\hat{x}^3_2)^\top=(1,0)^\top$ are two optimal solutions to problem $(P_{p_{sub1}})$.

By a similar discussion, we obtain that $\hat{x}^2_1=0$ and $\hat{x}^4_1=-1$ are two optimal solutions to problem $(\tilde{P}_{p_{sub2}})$ (see the black curve in Figure \ref{examples1-2}), which further implies that $\hat{\mathbf{x}}^2=(\hat{x}^2_1,\hat{x}^2_2)^\top=(0,1)^\top$ and $\hat{\mathbf{x}}^4=(\hat{x}^4_1,\hat{x}^4_2)^\top=(-1,0)^\top$ are two optimal solutions to problem $(P_{p_{sub2}})$.

Noticing that the objective values of sub-problem $(P_{p_{sub1}})$ at the points $\hat{\mathbf{x}}^1$ and $\hat{\mathbf{x}}^3$ as well as the objective values of sub-problem $(P_{p_{sub2}})$ at the points $\hat{\mathbf{x}}^2$ and $\hat{\mathbf{x}}^4$ are $1$, we obtain that the $\ell_p$-minimization (\ref{e-exam1-2}) has four optimal solutions: $\hat{\mathbf{x}}^1$, $\hat{\mathbf{x}}^2$, $\hat{\mathbf{x}}^3$, $\hat{\mathbf{x}}^4$ whenever $p\in (0,p^\ast)$.

From the above discussions, we can see that for any fixed $p\in (0,p^\ast)$, all optimal solutions of the $\ell_p$-minimization (\ref{e-exam1-2}) are also optimal solutions of the $\ell_0$-minimization (\ref{e-exam1-1}), which is consistent with our theoretical result.

\begin{Remark}
(i) The corresponding $\ell_1$-minimization of the $\ell_0$-minimization (\ref{e-exam1-1}) is
\begin{eqnarray}\label{e-exam1-3}
\min\limits_{x\in \mathbb{R}^2}\| x \|_1\quad \textrm{s.t.}\;\; |x_1-x_2|=1.
\end{eqnarray}
Obviously, the optimal solution set of this problem is $\{\textbf{x}\in \mathbb{R}^2| |x_1-x_2|=1, \|\textbf{x}\|_\infty\leqslant 1\}$, and hence, the $\ell_1$-minimization (\ref{e-exam1-3}) has infinite optimal solutions which are not optimal solutions to the $\ell_0$-minimization (\ref{e-exam1-1}). This implies that $p^\ast$ we found in Example \ref{example1-lp-solve-sparse} is the biggest one satisfying the desired property. (ii) From Example \ref{example1-lp-solve-sparse}, it is easy to see that for any fixed $p\in (0,p^\ast)$, optimal solutions of the $\ell_p$-minimization (\ref{e-exam1-2}) coincide with optimal solutions of the $\ell_0$-minimization (\ref{e-exam1-1}).
\end{Remark}

\begin{Example}\label{example2-lp-solve-sparse}
Consider a 2-dimensional problem:
\begin{eqnarray}\label{e-exam2-1}
\min\limits_{x\in \mathbb{R}^2}\| x \|_0\quad \textrm{s.t.}\;\; |5x_1+x_2|=2.
\end{eqnarray}
\end{Example}

\begin{figure}[ht]
\centering
\includegraphics[width=2.5in]{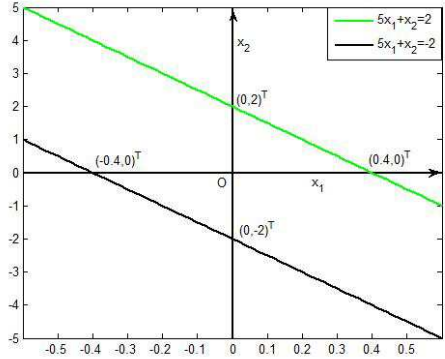}
\caption{The problem in Example \ref{example2-lp-solve-sparse} has four sparse solutions.}
\label{examples2-1}
\end{figure}

First, it is easy to show that this problem has four sparse solutions: $\hat{\mathbf{x}}^1:=(0.4,0)^\top$, $\hat{\mathbf{x}}^2:=(-0.4,0)^\top$, $\hat{\mathbf{x}}^3:=(0,2)^\top$ and $\hat{\mathbf{x}}^4:=(0,-2)^\top$ (See Figure \ref{examples2-1}).

Second, we compute the value of $p^\ast$. It is easy to see that the sparsity level of the $\ell_0$-minimization (\ref{e-exam2-1}) is $s=1$; $\Phi=(5,1)\in\mathbb{R}^{1\times2}$; $\epsilon\in\{-1,1\}$ and $\mathbf{b}=2$. Denote $\epsilon_1:=-1$ and $\epsilon_2:=1$, then  $\mathbf{b}_{\epsilon_{1}}=\mathbf{b}\circ\epsilon_{1}=-2$ and $\mathbf{b}_{\epsilon_{2}}=\mathbf{b}\circ\epsilon_{2}=2$. Similar to Example \ref{example1-lp-solve-sparse}, we can obtain that
\begin{eqnarray*}
c_0=\max\{c^{\epsilon_1}_0,c^{\epsilon_2}_0\}=\max\{2,2\}=2,\quad
c_1=\max\{c^{\epsilon_1}_p,c^{\epsilon_2}_p\}=\max\{\frac{10}{13},\frac{10}{13}\}=\frac{10}{13},
\end{eqnarray*}
and hence, $r_0=\max\{c_0,c_1\}=\max\{2,\frac{10}{13}\}=2$. It is easy to see that the set $\mathbb{T}$ has four pseudo-extreme points $(-0.4, 0, 0.4, 0)^{\top}$ , $(0,-2, 0, 2)^{\top}$, $(0.4, 0, 0.4, 0)^{\top}$ , $(0, 2, 0, 2)^{\top}$, and hence, $r_m=0.4$. Furthermore, $\frac{s+1}{s}=2<\frac{r_0}{r_m}=\frac{2}{0.4}=5$.
It follows from (\ref{p-star-formula}) that
$$p^\ast=\frac{\ln (s+1)-\ln s}{\ln r_0-\ln r_m }=\frac{\ln 2}{\ln 5}\approx 0.4307.$$

Third, for any fixed $p\in (0,p^\ast)$, we consider the corresponding $\ell_p$-minimization:
\begin{eqnarray}\label{e-exam2-2}
\min\limits_{x\in \mathbb{R}^2}\| x \|_p^p\quad \textrm{s.t.}\;\; |5x_1+x_2|=2,
\end{eqnarray}
and find all optimal solutions of this problem.

The $\ell_p$-minimization (\ref{e-exam2-2}) is equivalent to
$$
\min\{\underbrace{\min\limits_{5x_1+x_2=2}\| \mathbf{x} \|_p^p}\limits_{(P_{p_{sub1}})},\;\underbrace{\min\limits_{5x_1+x_2=-2}\| \mathbf{x } \|_p^p}\limits_{(P_{p_{sub2}})}\}.
$$
Both sub-problems can be equivalently transformed into one-variable unconstrained problems, i.e.,
\begin{eqnarray}\label{exam2-subp}
(\tilde{P}_{p_{sub1}}):\;\; \min\limits_{x_1\in\mathbb{R}}|x_1|^p+|2-5x_1|^p\quad
\mbox{\rm and}\quad
(\tilde{P}_{p_{sub2}}):\;\; \min\limits_{x_1\in\mathbb{R}}|x_1|^p+|-2-5x_1|^p.
\end{eqnarray}

\begin{figure}[ht]
\centering
\includegraphics[width=2.5in]{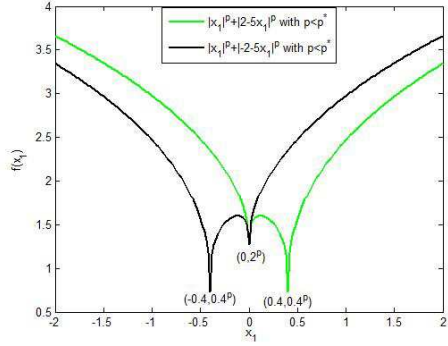}
\caption{Graphs of objective functions of two subproblems in (\ref{exam2-subp}) with one-variable $x_1$.}
\label{examples2-2}
\end{figure}

For sub-problem $(\tilde{P}_{p_{sub1}})$, we divide $\mathbb{R}$ into three parts: $\mathbb{R}=\mathbb{D}_1\bigcup\mathbb{D}_2\bigcup\mathbb{D}_3$
where $\mathbb{D}_1=(-\infty,0]$, $\mathbb{D}_2=[0,0.4]$ and $\mathbb{D}_3=[0.4,+\infty)$; and accordingly, we denote $f_1(x_1)=(-x_1)^p+(2-5x_1)^p,f_2(x_1)=x_1^p+(2-5x_1)^p$ and $f_3(x_1)=x_1^p+(5x_1-2)^p$. Then, problem $(\tilde{P}_{p_{sub1}})$ is equivalent to
$$
\min\{\underbrace{\min\limits_{x_1\in\mathbb{D}_1}f_1(x_1)}\limits_{(\tilde{P}_{p^1_{sub1}})},
\underbrace{\min\limits_{x_1\in\mathbb{D}_2}f_2(x_1)}\limits_{(\tilde{P}_{p^2_{sub1}})},
\underbrace{\min\limits_{x_1\in\mathbb{D}_3}f_3(x_1)}\limits_{(\tilde{P}_{p^3_{sub1}})}\}.
$$
We now consider the above three sub-problems, respectively.
\begin{itemize}
  \item [(i)] For problem $(\tilde{P}_{p^1_{sub1}})$. Since $f^\prime_1(x_1)=(-1)p(-x_1)^{p-1}+(-5)p(2-5x_1)^{p-1}<0$
      on $\mathbb{D}_1$, it follows that $f_1(\cdot)$ is a strictly decreasing function on $\mathbb{D}_1$, and hence,
      $\min\limits_{x_1\in\mathbb{D}_1}f_1(x_1)=f_1(0)=2^p.$
 \item [(ii)] For problem $(\tilde{P}_{p^2_{sub1}})$. It is easy to see that on $\mathbb{D}_2$,
     \begin{eqnarray*}
     \begin{array}{l}
     f^\prime_2(x_1)=px_1^{p-1}+(-5)p(2-5x_1)^{p-1}, \\ f^{\prime\prime}_2(x_1)=p(p-1)x_1^{p-2}+25p(p-1)(2-5x_1)^{p-2}<0,
     \end{array}
     \end{eqnarray*}
     which implies that the function $f_1(\cdot)$ is strictly concave on $\mathbb{D}_2$. Thus, the optimal value is obtained on the endpoints of $\mathbb{D}_2$, that is
     $$
     \min\limits_{x_1\in\mathbb{D}_2}f_2(x_1)=\min\{f_2(0),f_2(0.4)\}=f_2(0.4)=0.4^p.
     $$
\item [(iii)] For problem $(\tilde{P}_{p^3_{sub1}})$. Since
    $f^\prime_3(x_1)=px_1^{p-1}+5p(5x_1-2)^{p-1}>0$,
    it follows that $f_3(\cdot)$ is a strictly increasing function on $\mathbb{D}_3$, and hence,       $\min\limits_{x_1\in\mathbb{D}_3}f_3(x_1)=f_3(0.4)=0.4^p$.
\end{itemize}
Combining cases (i)-(iii), we obtain that $\hat{x}_1=0.4$ is the unique optimal solution to problem $(\tilde{P}_{p_{sub1}})$ (see the green curve in Figure \ref{examples2-2}). Furthermore,  $\hat{\mathbf{x}}^1=(\hat{x}_1,\hat{x}_2)^\top=(0.4,0)^\top$ is the unique optimal solution to problem $(P_{p_{sub1}})$.

By a similar discussion, we obtain that $\hat{x}_1=-0.4$ is the unique optimal solution to problem $(\tilde{P}_{p_{sub2}})$ (see the black curve in Figure \ref{examples2-2}), which further implies that $\hat{\mathbf{x}}^2=(\hat{x}_1,\hat{x}_2)^\top=(-0.4,0)^\top$ is the unique optimal solution to problem $(P_{p_{sub2}})$.

Noticing that both the objective value of problem $(P_{p_{sub1}})$ at the point $\hat{\mathbf{x}}^1$ and the objective value of problem $(P_{p_{sub2}})$ at the point $\hat{\mathbf{x}}^2$ are the same, i.e., $0.4^p$, we obtain that the $\ell_p$-minimization (\ref{e-exam2-2}) has two optimal solutions: $(0.4,0)^{\top}$ and  $(-0.4,0)^{\top}$.

From the above discussions, we can see that for any given $p\in (0,p^\ast)$, all optimal solutions of the $\ell_p$-minimization (\ref{e-exam2-2}) are also optimal solutions of the $\ell_0$-minimization (\ref{e-exam2-1}), which is consistent with our theoretical result.

\begin{Remark}
(i) From the above discussions, it is easy to see that the value of $p^\ast$ obtained in Example \ref{example2-lp-solve-sparse} is not the biggest one satisfying the desired property. In fact, the value of $p^\ast$ we obtained is $p^\ast=\frac{\ln 2}{\ln 5}\approx 0.4307$, however, it is not difficult to show that for any fixed $p\in (0,1)$, all optimal solutions of the $\ell_p$-minimization (\ref{e-exam2-2}) are optimal solutions of the $\ell_0$-minimization (\ref{e-exam2-1}). So, it is worthy of further investigation to improve the construction of $p^\ast$. (ii) The $\ell_0$-minimization (\ref{e-exam2-1}) has four optimal solutions and the $\ell_p$-minimization (\ref{e-exam2-2}) has two optimal solutions. We obtain that for any fixed $p\in (0,p^\ast)$ with $p^\ast=\frac{\ln 2}{\ln 5}$, all optimal solutions of the $\ell_p$-minimization (\ref{e-exam2-2}) solve the $\ell_0$-minimization (\ref{e-exam2-1}). It is easy to see that the $\ell_p$-minimization (\ref{e-exam2-2}) has also two local optimal solutions except two global optimal solutions mentioned above; and they are just other two optimal solutions of the $\ell_0$-minimization (\ref{e-exam2-1}).
\end{Remark}

\section{Conclusions}

In this paper,  we study the relationship between the $\ell_0$-minimization and the corresponding $\ell_p$-minimization for PCS in the real context.
We show that there exists a fixed constant $p^\ast>0$ such that every optimal solution to the $\ell_p$-minimization with any fixed $p\in (0, p^\ast)$ also solves the $\ell_0$-minimization.
Moreover, we show that such a constant $p^\ast$ can be expressed by the sparsity level of the concerned problem and the related given coefficient matrix and vector.  These provide a theoretical basis for solving PCS in the real context via solving the corresponding $\ell_p$-minimization. For the purpose of practical applications, effective algorithms need to be developed.

%
%

\end{document}